\documentclass{amsart}
\def \rr {\mathbb{R}}
\def \rn {\mathbb{R}^n}
\def \dunp {D_1^p(\mathbb{R}^n)}
\def \crit {p^\star}
\def \crits {p^\star\!(s)}
\def \nn {\mathbb{N}}
\newtheorem{thm}{Theorem}
\newtheorem{prop}{Proposition}
\newtheorem{claim}{Claim}[section]
\author{Roberta Filippucci}
\author{Patrizia Pucci}
\author{Fr\'ed\'eric Robert}
\address{R. Filippucci,
P. Pucci: Dipartimento di Matematica e Informatica,
Universit\`a degli Studi di Perugia, via Vanvitelli 1, 06123 Perugia, Italy }
\email{filippucci@dipmat.unipg.it,
pucci@dipmat.unipg.it}
\address{F. Robert: Laboratoire J.--A.Dieudonn\'e, Universit\'e de Nice--Sophia Antipolis, 06108 Nice Cedex 2, France}
\email{frobert@math.unice.fr}
\title[Multiple critical nonlinearities]{On a $p$--Laplace equation\\
with multiple critical nonlinearities}
\date{September 18, 2008.}

\medskip

\begin{document}

\begin{abstract}
Using the Mountain--Pass Theorem of Ambrosetti and Rabinowitz we
prove that $-\Delta_p
u-\mu|x|^{-p}{u^{p-1}}=|x|^{-s}{u^{\crits-1}}+u^{\crit-1}$ admits a
positive weak solution in $\rn$ of class $\dunp\cap
C^1(\rn\setminus\{0\})$, whenever $\mu<\mu_1$, and
$\mu_1=[(n-p)/p]^p$. The technique is based
on the existence of extremals of some Hardy--Sobolev type embeddings
of independent interest. We also show that if $u\in\dunp$ is a weak
solution in $\rn$ of $-\Delta_p
u-\mu|x|^{-p}{|u|^{p-2}u}=|x|^{-s}{|u|^{\crits-2}u}+|u|^{q-2}u$,
then $u\equiv0$ when either $1<q<\crit$, or $q>\crit$ and $u$ is
also of class $L^\infty_\text{\scriptsize{loc}}(\rn\setminus\{0\})$.
\end{abstract}

\maketitle

\tableofcontents

\section{Introduction}

In this paper,  we are interested in weak solutions $u\in\dunp$,
$u\geq 0$ a.e, of the {\it double} critical equation of Emden-Fowler
type
\begin{equation}\label{eq:th}
-\Delta_p
u-\mu\frac{u^{p-1}}{|x|^p}=u^{\crit-1}+\frac{u^{\crits-1}}{|x|^s}\quad\mbox{
in }\rn,
\end{equation}
where $\Delta_p:=\hbox{div}(|\nabla u|^{p-2}\nabla u)$ is the
$p$--Laplace operator, $n\geq 2$ is an integer, $\mu$ is a real
parameter, $p\in (1,n)$ and $\crit:={np}/(n-p)$, while $s\in (0,p)$
and  $\crits:={p(n-s)}/(n-p)$. The space $\dunp$ is defined as the
completion of $C^\infty_c(\rn)$, the set of smooth compactly
supported function on $\rn$, for the norm
$$u\mapsto\Vert\nabla u\Vert_p,$$
where here and in the sequel, $\Vert\cdot\Vert_q$ denotes the
$L^q$--norm on the Lebesgue space $L^q(\rn)$.

Throughout the paper, we say that $u\in \dunp$ is a {\it weak
solution} of $-\Delta_pu=f$, where $f\in L^1_{\text{\scriptsize{\rm
loc}}}(\rn)$, if
$$\int_{\rn}|\nabla u|^{p-2}(\nabla u,\nabla\varphi)\, dx=\int_{\rn} f\varphi\, dx$$
for all $\varphi\in C_c^\infty(\rn)$.
\smallskip

Existence and non--existence, as well as qualitative properties, of
non--trivial non--negative solutions for elliptic equations with
singular potentials were recently studied by several authors, but,
essentially, only with a solely critical exponent. We refer,
e.g., in bounded domains and for $p=2$ to \cite{ch, gr1, gr2, kp0, kp},
and for general $p>1$ to \cite{ck,dh,gy,h}; while in $\rn$ and for
$p=2$ to \cite{cw, ft, li, terr}, and for general $p>1$ to \cite{afp, gm,
patraf}, and the references therein. The large literature on
$p$--Laplacian equations in the entire $\rn$ differs somehow for the
nonlinear structure, objectives and methods  from those presented in
this paper.

\medskip\noindent Indeed, the combination of the two critical exponents induces more
subtilizes and difficulties. When just one critical exponent is
involved, there are solutions to the corresponding equations (see
for instance \cite{patraf}): in general, these solutions are
radially symmetrical with respect to a point of the domain ($0$ in
general) and are explicit. In our context, very few is known: yet, we refer to an interessant approach by Kang and Li \cite{kl}.

\medskip\noindent A natural strategy is to construct the  solutions
of \eqref{eq:th} as critical points of a suitable functional via the
mountain--pass lemma of Ambrosetti and Rabinowitz. Due to the
invariance of \eqref{eq:th} by the conformal one parameter
transformation group
\begin{equation}\label{eq:transfo}
\left\{\begin{array}{llll}
T_r: & \dunp & \to & \dunp\\
     & u & \mapsto & \left[x\mapsto r^{(n-p)/p}u(rx)\right]
\end{array}\right\},\; r>0,
\end{equation}
it is well--known that the mountain--pass lemma does not yield
critical points, but only Palais--Smale sequences. The main issue of
the paper is to understand the behavior of these Palais--Smale
sequences. Indeed, the principal difficulty here is that there is an
asymptotic competition between the energy carried by the two
critical nonlinearities. If one dominates the other, then there is
vanishing of the weakest one and one recovers solutions to an
equation with only one critical nonlinearity: in this situation, we
do not get solutions of equation \eqref{eq:th}. Therefore, the
crucial point here is to avoid the domination of one term on the
other.

\medskip\noindent Sections~\ref{sec:ps}, \ref{sec:struc} and~\ref{sec:proof} of the paper are devoted
to the proof of the following main existence result:
\begin{thm}\label{th:main} For any $\mu\in(-\infty,\mu_1)$, $\mu_1:=[(n-p)/p]^p$, and $s\in (0,p)$,
there exists a positive weak solution of \eqref{eq:th}. More
precisely, there exists $u\in D_1^p(\rn)\cap C^1(\rn\setminus\{0\})$
such that $u>0$ in $\rn\setminus\{0\}$ and $u$ solves \eqref{eq:th}
weakly in $\rn.$
\end{thm}

\medskip\noindent Theorem~\ref{th:main} is proved via the choice of a suitable energy level
for the mountain--pass lemma: with this choice, a careful analysis
of concentration  allows us to show that there is a balance between
the energies of the two nonlinearities mentioned above, and
therefore none can dominate the other. There we make a full use of
the conformal invariance of \eqref{eq:th} under the transformation
\eqref{eq:transfo}; this guarantees the convergence to a solution to
\eqref{eq:th}. As an offshoot of this analysis, we prove that the
blow--up energy is quantized for both nonlinearities.\par

\medskip\noindent The choice of the energy level involves the best
constants in the Hardy--Sobolev inequalities (see \eqref{def:K} and
\eqref{def:Ks} of Section~\ref{sec:ps}). We are then led to
considering the possible extremals for them. As far as we know, the
result in its full generality, that we need, does not appear in the
literature: therefore, for the sake of completeness, we prove the
existence of extremals when $s>0$ in the Appendix given in
Section~\ref{sec:app}. Concerning the case $s= 0$, there is no
extremal in general when $\mu<0$ and the analysis relies on the radial case and is made in the Appendix given in Section~\ref{sec:neg:mu}. For details
concerning the extremals in the case $s=0$ we remind to both Sections~\ref{sec:app}
and~\ref{sec:neg:mu}.

\medskip\noindent It is to be noticed that the exponents $\crit$ and $\crits$ are exactly the ones
that make the equation invariant under the  transformation group
\eqref{eq:transfo}. One can therefore naturally wonder what happens
for different exponents:  in  Section \ref{sec:poho}, we present a
non--existence theorem, when $q\neq \crit$, cf.
Theorem~\ref{th:nonex} and Claims~\ref{claimp4}--\ref{claimp5} (we
also refer to \cite{patraf} for other nonexistence results in the
same spirit). In particular, in general, there is no solution to the
corresponding equation (except the null one) when one takes
exponents different from $\crit$ and $\crits$ in \eqref{eq:th}.

\medskip\noindent The paper is organized as follows: in Sections~\ref{sec:ps}, \ref{sec:struc} and~\ref{sec:proof}
we prove Theorem~\ref{th:main} when $\mu\geq 0$. In Section~\ref{sec:poho} we deal
with  the non--existence result in the spirit of Pohozaev. In
Section~\ref{sec:app}, we prove the existence of extremals for some
Hardy--Sobolev type embeddings, see Theorem~\ref{th:ext}. While
Section~\ref{sec:neg:mu} deals with the situation in which $\mu<0$.

\section{Preliminaries and construction of the appropriate Palais--Smale sequence}\label{sec:ps}
Clearly equation \eqref{eq:th} is related to some specific
functional embeddings and inequalities. The standard Hardy
inequality asserts that $\dunp$ is embedded in the weighted space
$L^p(\rn,|x|^{-p})$ and that this embedding is continuous: more
precisely,
\begin{equation}\label{ineq:hardy}
\mu_1\int_{\rn}\frac{|u|^p}{|x|^p}\, dx\leq \int_{\rn}|\nabla
u|^pdx,\qquad \mu_1:=\left(\frac{n-p}{p}\right)^p,
\end{equation}
for all $u\in \dunp$.  Moreover, the constant $\mu_1$ is optimal. If
$\mu<\mu_1$, it follows from the Hardy inequality
$\eqref{ineq:hardy}$ that
$$\Vert u\Vert:=\left(\int_{\rn}|\nabla u|^p\, dx-\mu\int_{\rn}\frac{|u|^p}{|x|^p}\, dx\right)^{1/p}$$
is well--defined on $\dunp$. Note that, $\Vert\cdot\Vert$ is {\it
comparable} to the norm $\Vert\nabla\cdot\Vert_p$ since the
following inequalities hold
\begin{equation}\label{comp:norm}
\left(1-\frac{\mu_+}{\mu_1}\right)\Vert\nabla u\Vert_p^p\leq \Vert
u\Vert^p\leq \left(1+\frac{\mu_-}{\mu_1}\right)\Vert\nabla
u\Vert_p^p
\end{equation}
for any $u\in \dunp$, where $\mu_+=\max\{\mu,\,0\}$ and
$\mu_-=\max\{-\mu,\,0\}$.
\smallskip

It follows from Sobolev's embedding theorem that $\dunp$ is
continuously embedded  in $L^{\crit}(\rn)$ where
$\crit:={np}/(n-p)$. Therefore, there exists $C>0$ such that $\Vert
u\Vert_{\crit}\leq C\Vert u\Vert$. Taking $C$ as small as possible,
we define the optimal constant $K(n,p,\mu,0)>0$ associated to this
embedding as

\begin{equation}\label{def:K}
\frac{1}{K(n,p,\mu,0)}:=\inf_{u\in \dunp\setminus\{0\}}
\frac{\int_{\rn}|\nabla u|^p\, dx-\mu\int_{\rn}{|u|^p}|x|^{-p}\, dx}
{\left(\int_{\rn}|u|^{\crit}\, dx\right)^{{p}/{\crit}}},
\end{equation}
that is $C^p=K(n,p,\mu,0)$. Combining the Hardy inequality and the
Sobolev inequality, we obtain the Hardy--Sobolev inequality. Indeed,
let $s\in (0,p)$ be a real number: then $\dunp$ is continuously
embedded in the weighted space $L^{\crits}(\rn,|x|^{-s})$, where
$\crits:={p(n-s)}/(n-p)$. Here again, taking the smallest constant
associated to this embedding, we let

\begin{equation}\label{def:Ks}
\frac{1}{K(n,p,\mu,s)}:=\inf_{u\in
\dunp\setminus\{0\}}\frac{\int_{\rn}|\nabla u|^p\, dx
-\mu\int_{\rn}{|u|^p}|x|^{-p}\,
dx}{\left(\int_{\rn}{|u|^{\crits}}|x|^{-s}\,
dx\right)^{{p}/{\crits}}}.
\end{equation}

Let the functional $\Phi$ defined on $\dunp$ as follows:
$$\Phi(u):=\frac{1}{p}\Vert u\Vert^p-\frac{1}{\crit}\int_{\rn}(u_+)^{\crit}\, dx-\frac{1}{\crits}\int_{\rn}\frac{(u_+)^{\crits}}{|x|^s}\, dx$$
for $u\in\dunp$. Here and in the sequel, $u_+=\max\{u,0\}$. It
follows from the Hardy, Sobolev and Hardy--Sobolev embeddings that
$\Phi$ is well--defined and that $\Phi\in C^1(\dunp)$. Note that a
positive weak solution to \eqref{eq:th} is a nontrivial critical
point of $\Phi$; and we actually show, in the proof of
 Claim~\ref{claim33}, that a nonnegative nontrivial weak
 limit of a Palais--Smale sequence of $\Phi$ is a positive solution
 of \eqref{eq:th} by the Tolksdorf regularity theory
 \cite{tolk} and the Vazquez strong maximum principle~\cite{vaz}.

In this section, we prove the following:

\begin{prop}\label{prop:ps} Assume
\begin{equation}\label{or}
\mu\in[0,\mu_1)\hbox{ and }s\in[0,p).
\end{equation}
Then
there exists $(u_k)_{k\in\nn}\in\dunp$ such that
$$\lim_{k\to \infty}\Phi'(u_k)= 0\quad\hbox{strongly in }(\dunp)'\quad\hbox{and}\quad\lim_{k\to \infty}\Phi(u_k)=c,$$
where
\begin{equation}\label{bnd:c}0<c<c_\star:=
\min\left\{\frac1n K(n,p,\mu,0)^{-n/p},
\frac{p-s}{p(n-s)}K(n,p,\mu,s)^{-(n-s)/(p-s)}\right\}.
\end{equation}
\end{prop}

\medskip\noindent Note that ${1}/{p}- {1}/{\crit}=1/n$, ${\crit}/(\crit-p)=n/p$,
${1}/{p}- {1}/{\crits}=(p-s)/p(n-s)$ and
${\crits}/(\crits-p)=(n-s)/(p-s)$. The proof of Proposition
\ref{prop:ps} uses the following version of the Mountain--Pass
lemma:
\begin{thm}[Ambrosetti and Rabinowitz, \cite{ar}]\label{th:mpl}
Let $(V,N)$ be a Banach space and let $F\in C^1(V)$. We assume
that\par (i) $F(0)=0$,\par (ii) There exist $\lambda,R>0$ such that
$F(u)\geq\lambda$ for all $u\in V$, with $N(u)=R$,\par (iii) There
exists $v_0\in V$ such that $\limsup_{t\to \infty}F(tv_0)<0$.\par
\noindent Let $t_0>0$ be such that $N(t_0v_0)>R$ and $F(t_0v_0)<0$
and let
$$c:=\inf_{\gamma\in\Gamma}\,\sup_{t\in [0,1]}F(\gamma(t)),$$
where
$$\Gamma:=\{\gamma\in C^0([0,1],V)\,/\, \gamma(0)=0\hbox{ and }\gamma(1)=t_0v_0\}.$$
Then there exists a Palais--Smale sequence at level $c$, that is
there exists a sequence $(u_k)_{k\in\nn}\in V$ such that
$$\lim_{k\to \infty}F(u_k)=c\quad\hbox{ and }\quad\lim_{k\to \infty}F'(u_k)= 0\quad\hbox{strongly in }V'.$$
\end{thm}

\begin{claim}\label{claim1}
The functional $\Phi$ verifies the hypotheses of the Mountain--Pass
lemma at any $u\in\dunp$, with $u_+\not\equiv 0$.\end{claim}

\medskip\noindent{\it Proof of Claim~$\ref{claim1}$:} Clearly $\Phi\in C^1(\dunp)$ and $\Phi(0)=0$.
Using the definition of the best constants in \eqref{def:K},
\eqref{def:Ks}, we get that
$$\begin{aligned}
\Phi(u)&\geq \frac{1}{p}\Vert u\Vert^p-
\frac{K(n,p,\mu,0)^{{\crit}/{p}}}{\crit}\Vert u\Vert^{{\crit}}
-\frac{K(n,p,\mu,s)^{{\crits}/{p}}}{\crits}\Vert u\Vert^{{\crits}}\\
&=  \left(\frac{1}{p}-\frac{K(n,p,\mu,0)^{{\crit}/{p}}}{\crit}\Vert
u\Vert^{{\crit}-p} -\frac{K(n,p,\mu,s)^{{\crits}/{p}}}{\crits}\Vert
u\Vert^{{\crits}-p}\right) \Vert u\Vert^p.
\end{aligned}$$
Then, since \eqref{comp:norm} holds and since $p<\crits\le\crit$ being
$s\in[0,p)$, there exists $R>0$ such that $\Phi(u)\ge\lambda$ for
all $u\in\dunp$ such that $\Vert\nabla u\Vert_p=R$: point (ii) of
Theorem~\ref{th:mpl} is satisfied. Moreover, given any $u\in\dunp$, with $u_+\not\equiv 0$, we have that
$$\lim_{t\to \infty}\Phi(t u)=-\infty.$$
We then, let $t_u>0$ be such that $\Phi(t
 u)<0$ for $t\geq t_u$ and $\Vert\nabla( t_u u)\Vert_p>R$. Consider
$$\Gamma_u:=\{\gamma\in C^0([0,1],\dunp)\,/\, \gamma(0)=0\hbox{ and }\gamma(1)=t_u u\}$$
and
$$c_u:=\inf_{\gamma\in\Gamma_u}\,\sup_{t\in [0,1]}\Phi(\gamma(t)).$$
Then the hypotheses of Theorem~\ref{th:mpl} are satisfied. This ends
the proof of Claim~\ref{claim1}.\hfill$\Box$\par
\medskip\noindent It follows from Theorem~\ref{th:mpl} that there exists $(u_k)_{k\in\nn}\in\dunp$ such that
$$\lim_{k\to \infty}\Phi(u_k)=c_u\quad\hbox{ and }\quad\lim_{k\to \infty}\Phi'(u_k)= 0\quad\hbox{strongly in }(\dunp)'.$$
Moreover,  from the definition of $c_u$ it is also clear that
$c_u\geq\lambda$, and so
$$c_u>0$$ for all $u\in\dunp\setminus\{0\}$.

\begin{claim}\label{claim2} Assume  \eqref{or}. Then there exists $u\in\dunp\setminus\{0\}$ such that $u\geq 0$ and
\begin{equation}\label{ineq:step1}c_u<\frac1nK(n,p,\mu,0)^{-n/p}.\end{equation}\end{claim}

\medskip\noindent{\it Proof of Claim~$\ref{claim2}$:} By \eqref{or}, let $u\in\dunp\setminus\{0\}$ be a non-negative extremal for $1/K(n,p,\mu,0)$ in
\eqref{def:K} $($see Theorem~$\ref{th:ext}$ in
Section~$\ref{sec:app})$. Since $u=u_+$, by the definition of $t_u$
and the fact that $c_u>0$, we have
$$c_u\leq \sup_{t\geq 0}\Phi(t u)\leq \sup_{t\geq 0}f(t),$$
where
$$f(t):=\frac{t^p}{p}\Vert u\Vert^p-\frac{t^{\crit}}{\crit}\int_{\rn}|u|^{\crit}\, dx$$
for all $t\geq 0$. Straightforward computations yield
$$c_u\leq \left(\frac{1}{p}-\frac{1}{\crit}\right)\left(\frac{\Vert u\Vert^p}{\left(\int_{\rn}|u|^{\crit}\, dx\right)^{{p}/{\crit}}}\right)^{{\crit}/(\crit-p)}
=\frac{1}{n}K(n,p,\mu,0)^{-n/p},
$$
since $u$ is a non-negative extremal for \eqref{def:K}. Hence, if
equality would hold in \eqref{ineq:step1}, then $0<c_u=\sup_{t\geq
0}\Phi(t u)= \sup_{t\geq 0}f(t)$. Letting $t_1,t_2>0$ be points
where the two suprema are attained respectively, we get that
$$f(t_1)-\frac{t_1^{\crits}}{\crits}\int_{\rn}\frac{|u|^{\crits}}{|x|^s}\, dx=f(t_2),$$
that is $f(t_2)< f(t_1)$, being $u_+\not\equiv0$ and $t_1>0$. This
gives the required contradiction and the claim is proved when \eqref{or} holds.\hfill$\Box$

\begin{claim}\label{claim3} Assume  \eqref{or}. There exists $u\in\dunp\setminus\{0\}$ such that $u\ge0$ and
$$0<c_u<c_\star,$$
where $c_\star$ is defined in \eqref{bnd:c}.\end{claim}

\noindent{\it Proof of Claim~$\ref{claim3}$:} In case
$$\frac1n K(n,p,\mu,0)^{-n/p}\le
\frac{p-s}{p(n-s)}K(n,p,\mu,s)^{-(n-s)/(p-s)},$$ we take
$u\in\dunp\setminus\{0\}$ as in Claim~\ref{claim2} to get the
result. Otherwise we take
$u\in\dunp\setminus\{0\}$ a non--negative extremal for
\eqref{def:Ks} $($which exists by Theorem~$\ref{th:ext}$ of
Section~$\ref{sec:app})$ and proceed as in the first part of the
proof of Claim~\ref{claim2}, with $f$ replaced by
$$\tilde f(t):=\frac{t^p}{p}\Vert u\Vert^p-\frac{t^{\crits}}{\crits}\int_{\rn}\frac{|u|^{\crits}}{|x|^s}\, dx,$$
which gives now the contradiction
$$\tilde f(t_1)-\frac{t_1^{\crit}}{\crit}\int_{\rn}|u|^{\crit} dx=\tilde f(t_2).$$
This proves Claim~\ref{claim3}. \hfill$\Box$

\medskip\noindent Proposition \ref{prop:ps} is a consequence of Claims \ref{claim1} and \ref{claim3} for a suitable $u$ in $\dunp$.

\section{The structure of Palais--Smale sequence going to zero weakly}\label{sec:struc}
From now on, we assume that $s\in (0,p)$. We prove the following proposition:
\begin{prop}\label{prop:asymp} Let $(u_k)_{k\in\nn}\in\dunp$ be a Palais--Smale sequence at level $c\in (0, c_\star)$ as
in Proposition~$\ref{prop:ps}$, with $s\not=0$ in \eqref{or}. If $u_k\rightharpoonup 0$ weakly in
$\dunp$ as $k\to \infty$, then there exists
$\epsilon_0=\epsilon_0(n,p,\mu,s,c)>0$ such that
$$\hbox{either}\quad\lim_{k\to \infty}\int_{B_\delta(0)}(u_k)_+^{\crit}\, dx=0\quad\mbox{or}\quad\limsup_{k\to \infty}\int_{B_\delta(0)}(u_k)_+^{\crit}\, dx\geq\epsilon_0$$
for all $\delta>0$.
\end{prop}
\noindent The proof of Proposition \ref{prop:asymp} goes through
four claims.

\begin{claim}\label{claim21} Let $(u_k)_{k\in\nn}\in\dunp$ be a Palais--Smale sequence as in
Proposition~$\ref{prop:asymp}$. If $u_k\rightharpoonup 0$ weakly in
$\dunp$ as $k\to \infty$, then for all $\omega\subset\subset
\rn\setminus\{0\}$, up to a subsequence, we have that
\begin{eqnarray}
&&\lim_{k\to \infty}\int_\omega\frac{|u_k|^p}{|x|^p}\, dx=
\lim_{k\to \infty}\int_\omega\frac{|u_k|^{\crits}}{|x|^s}\, dx=0,\label{1}\\
&&\lim_{k\to \infty}\int_\omega |u_k|^{\crit}\, dx=\lim_{k\to
\infty}\int_\omega |\nabla u_k|^p\, dx=0.\label{2}
\end{eqnarray}\end{claim}

\noindent{\it Proof of Claim~$\ref{claim21}$:} Fix
$\omega\subset\subset\rn\setminus\{0\}$. Clearly the embedding
$\dunp\hookrightarrow L^q(\omega)$ is compact for $1\leq q<\crit$
and $|x|+|x|^{-1}$ is bounded on $\omega$. Hence \eqref{1} follows
at once, being $p<\crit$ and $\crits<\crit$ since $s\in(0,p)$ by
assumption.\par

\medskip\noindent Concerning the two equalities in \eqref{2},
let $\eta\in C^\infty_c(\rn\setminus\{0\})$ such that
$0\leq \eta\leq 1$ and $\eta_{|\omega}\equiv 1$. Since $\eta^pu_k\in\dunp$ for all $k\in\nn$, we get that
\begin{equation}\label{3}\langle\Phi'(u_k), \eta^p u_k\rangle=o(\Vert \eta^pu_k\Vert)=o(\Vert
u_k\Vert)=o(1)\end{equation} as $k\to \infty$, being
$(\Vert u_k\Vert)_{k\in\nn}$ bounded by the weak convergence of $(u_k)_{k\in\nn}$ in $\dunp$ and
\eqref{comp:norm}. Since $\lim_{k\to
\infty}\Vert u_k \Vert_{L^p(\text{\scriptsize{\rm
Supp}}|\nabla\eta|)}=0$ by \eqref{1} and $(\Vert\nabla u_k\Vert_p)_{k\in\nn}$ is bounded, we have as $k\to\infty$
$$\int_{\rn}|\nabla u_k|^{p-1}|\nabla(\eta^p)|\cdot |u_k|\, dx\le
\Vert \nabla u_k\Vert^{p-1}_p \Vert
u_k\Vert_{L^p(\text{\scriptsize{\rm Supp}}|\nabla\eta|)}=o(1),$$ and
so by~\eqref{3}
\begin{equation}\label{ineq:quanti}\begin{aligned}
o(1)&= \langle\Phi'(u_k), \eta^p u_k\rangle\\
&= \int_{\rn}|\nabla u_k|^{p-2}(\nabla u_k,\nabla(\eta^pu_k))\, dx-\int_{\rn}\eta^p(u_k)_+^{\crit}\, dx+o(1)\\
&=\int_{\rn}|\eta\nabla u_k|^{p}\, dx-\int_{\rn}\eta^p(u_k)_+^{\crit}\, dx \\
&\qquad
+O\left(\int_{\rn}|\nabla u_k|^{p-1}|\nabla(\eta^p)|\cdot |u_k|\, dx\right)+o(1) \\
&=\int_{\rn}|\nabla(\eta
u_k)|^{p}dx-\int_{\rn}\eta^p(u_k)_+^{\crit}\, dx
+o(1) \\
&\ge\Vert\eta u_ k\Vert^p-\int_{\rn}\eta^p(u_k)_+^{\crit}\, dx
+o(1),
\end{aligned}\end{equation}
since
\begin{equation}\label{io}\int_{\rn}|\nabla(\eta u_k)|^{p}dx=\int_{\rn}|\eta\nabla u_k|^{p}dx+o(1).\end{equation}
We prove \eqref{io}. Indeed, by the elementary inequality
$|\,|X+Y|^p-|X|^p|\le C_p(|X|^{p-1}+|Y|^{p-1})|Y|$ for all $X$,
$Y\in\rn$, we have $|\,|\nabla(\eta u_k)|^p-|\eta\nabla u_k|^p|\le
C_p(|\eta\nabla u_k|^{p-1} +|u_k\nabla\eta|^{p-1})|u_k\nabla\eta|$,
and by H\"older's inequality
$$\int_{\rn}|\eta\nabla u_k|^{p-1}|u_k\nabla\eta|\,dx\le
\Vert \nabla u_k\Vert^{p-1}_p \Vert
u_k\Vert_{L^p(\text{\scriptsize{\rm Supp}}|\nabla\eta|)}=o(1)$$ by
\eqref{1}, as well as $\int_{\rn}|u_k\nabla\eta|^{p}dx\le\Vert
u_k\Vert^p_{L^p(\text{\scriptsize{\rm Supp}}|\nabla\eta|)}=o(1)$.
This proves \eqref{io}.

\medskip\noindent Formula \eqref{ineq:quanti} above shows that
$$\Vert\eta u_ k\Vert^p\le\int_{\rn}(u_k)_+^{\crit-p}|\eta u_k|^p\, dx+o(1)$$
as $k\to \infty$. By H\"older's inequality and \eqref{def:K}, we
then have
$$\begin{aligned}
\Vert\eta  u_ k\Vert^p&\leq \left(\int_{\rn}(u_k)_+^{\crit}\,
dx\right)^{(\crit-p)/{\crit}}
\left(\int_{\rn}|\eta u_k|^{\crit}\, dx\right)^{ {p}/{\crit}}+o(1)\\
&\leq \left(\int_{\rn}(u_k)_+^{\crit}\,
dx\right)^{(\crit-p)/{\crit}}K(n,p,\mu,0)\, \Vert\eta u_ k\Vert^p
+o(1),
\end{aligned}$$
which gives
\begin{equation}\label{eq:1}
\left(1-\left(\int_{\rn}(u_k)_+^{\crit}\,
dx\right)^{(\crit-p)/{\crit}}K(n,p,\mu,0)\right) \Vert\eta  u_
k\Vert^p \leq o(1).
\end{equation}
Independently,
$\Phi(u_k)-\frac{1}{p}\langle\Phi'(u_k),u_k\rangle=c+o(\Vert
u_k\Vert)=c+o(1)$ as $k\to \infty$ since $(\Vert
u_k\Vert)_{k\in\nn}$ in bounded, which yields
\begin{equation}\label{balance:c}
\left(\frac{1}{p}-\frac{1}{\crit}\right)\int_{\rn}(u_k)_+^{\crit}\,
dx+\left(\frac{1}{p}-\frac{1}{\crits}\right)\int_{\rn}\frac{(u_k)_+^{\crits}}{|x|^s}\,
dx=c+o(1)
\end{equation}
as $k\to \infty$. Therefore,
\begin{equation}\label{upp:bnd:alpha}
\int_{\rn}(u_k)_+^{\crit}\, dx\leq {c}\,n+o(1)
\end{equation}
as $k\to \infty$. Plugging  \eqref{upp:bnd:alpha} into \eqref{eq:1}
we get that
$$\left(1-\left( c\,n\right)^{p/n}K(n,p,\mu,0)\right)
\Vert\eta u_ k\Vert^p\leq o(1)$$ as $k\to \infty$. The upper bound
\eqref{bnd:c} on $c$ yields
$$\lim_{k\to \infty}\Vert\eta u_ k\Vert^p =0,$$
and in turn by \eqref{def:K}
$$\lim_{k\to \infty}\int_{\rn}|\eta u_k|^{\crit}\, dx=0.$$
Since $\eta_{|\omega}\equiv 1$, these two latest inequalities and
\eqref{comp:norm} yield \eqref{2}. This proves
Claim~\ref{claim21}.\hfill$\Box$

\bigskip\noindent For $\delta>0$, we define
\begin{equation}\label{def:abc}\begin{gathered}
\alpha:=\limsup_{k\to \infty}\int_{B_\delta(0)}(u_k)_+^{\crit}\,dx;
\qquad \beta:=\limsup_{k\to \infty}\int_{B_\delta(0)}\frac{(u_k)_+^{\crits}}{|x|^s}\,dx\,;\\
\gamma:=\limsup_{k\to \infty}\int_{B_\delta(0)} \left(|\nabla
u_k|^{p}-\mu\frac{|u_k|^p}{|x|^p}\right)dx.
\end{gathered}\end{equation}
It follows from Claim~\ref{claim21} that these three quantities are
well--defined and independent of the choice of $\delta>0$.

\begin{claim}\label{claim22} Let $(u_k)_{k\in\nn}\in\dunp$ be a Palais--Smale sequence as in
Proposition~$\ref{prop:asymp}$, and let $\alpha$, $\beta$ and $\gamma$
be defined as in \eqref{def:abc}. If $u_k\rightharpoonup 0$ weakly
in $\dunp$ as $k\to \infty$, then
\begin{equation}\label{ineq:2.2}
\alpha^{ {p}/{\crit}}\leq K(n,p,\mu,0)\gamma\quad\hbox{ and
}\quad\beta^{ {p}/{\crits}}\leq K(n,p,\mu,s)\gamma.
\end{equation}\end{claim}

\noindent{\it Proof of Claim~$\ref{claim22}$:} Let $\eta\in
C^\infty_c(\rn)$ be such that $\eta_{|B_\delta(0)}\equiv 1$, with
$\delta>0$. Inequality \eqref{def:K} and Claim~\ref{claim21} yield
\begin{eqnarray*}
\left(\int_{\rn}|(\eta u_k)_+|^{\crit}\,
dx\right)^{{p}/{\crit}}&\leq &K(n,p,\mu,0)\,
\Vert\eta u_k\Vert^p\\
\left(\int_{B_\delta(0)}(u_k)_+^{\crit}\,
dx\right)^{{p}/{\crit}}&\leq &K(n,p,\mu,0)
\int_{B_\delta(0)}\left(|\nabla
u_k|^{p}-\mu\frac{|u_k|^p}{|x|^p}\right)\, dx+o(1)
\end{eqnarray*}
as $k\to \infty$. Letting $k\to \infty$, we get that $\alpha^{
{p}/{\crit}}\leq K(n,p,\mu,0)\gamma$. Similarly, we obtain the
second inequality of \eqref{ineq:2.2}. This proves Claim
\ref{claim22}.\hfill$\Box$

\begin{claim}\label{claim23} Let $(u_k)_{k\in\nn}\in\dunp$ be a Palais--Smale sequence as
in Proposition~$\ref{prop:asymp}$, and let $\alpha$, $\beta$ and
$\gamma$ be defined as in \eqref{def:abc}. If $u_k\rightharpoonup 0$
weakly in $\dunp$ as $k\to \infty$, then
$\gamma\leq\alpha+\beta$.\end{claim}

\medskip\noindent{\it Proof of Claim~$\ref{claim23}$:}
Let $\eta\in C^\infty_c(\rn)$ be such that
$\eta_{|B_\delta(0)}\equiv 1$. Since $\eta u_k\in \dunp$ and since
$\lim_{k\to \infty}\langle\Phi'(u_k),\eta u_k\rangle=0$, using
Claim~\ref{claim21} and the definitions of $\alpha$, $\beta$ and
$\gamma$ in \eqref{def:abc}, we get that $\gamma\leq\alpha+\beta$.
This proves Claim~\ref{claim23}.\hfill$\Box$

\medskip\noindent{\bf Proof of Proposition \ref{prop:asymp}}: Let $(u_k)_{k\in\nn}$ be as in
Proposition~\ref{prop:ps}, with $s\not=0$.
Claims~\ref{claim22} and~\ref{claim23} yield
\begin{equation}
\begin{gathered}\alpha^{ {p}/{\crit}}\leq K(n,p,\mu,0)\alpha +K(n,p,\mu,0)\beta,\\
\alpha^{ {p}/{\crit}}\left(1-K(n,p,\mu,0)\alpha^{ (\crit-p)/{\crit}}
\right)\leq K(n,p,\mu,0)\beta.\label{rel:alpha:beta}
 \end{gathered}
\end{equation}
Moreover, by \eqref{upp:bnd:alpha}, we obtain
\begin{equation}\label{upp:alpha}
\alpha\leq c\,n.
\end{equation}
Plugging \eqref{upp:alpha} into \eqref{rel:alpha:beta}, we have
$$\left(1-\left(c\,n\right)^{p/n}K(n,p,\mu,0)\right)
\alpha^{ {p}/{\crit}}\leq K(n,p,\mu,0)\beta.$$ By the upper bound
\eqref{bnd:c} on $c$  there exists $\delta_1$, depending on $n$,
$p$, $\mu$ and $c$, such that $\alpha^{ {p}/{\crit}}\leq
\delta_1\beta$. Similarly, there exists $\delta_2$, depending on
$n$, $p$, $\mu$, $s$ and $c$, such that $\beta^{ {p}/{\crits}}\leq
\delta_2\alpha$. In particular, it follows from these two latest
inequalities that there exists
$\epsilon_0=\epsilon_0(n,p,\mu,s,c)>0$ such that
\begin{equation}\label{4}
\hbox{either}\quad\alpha=\beta=0\quad\mbox{ or }\quad\{\alpha\geq
\epsilon_0\,\mbox{ and }\,\beta\geq \epsilon_0\}.
\end{equation}
By the definitions of $\alpha$ and $\beta$ given in \eqref{def:abc},
this proves Proposition~\ref{prop:asymp}. \hfill$\Box$

\section{Proof of Theorem \ref{th:main} in the case $\mu\ge0$}\label{sec:proof}
The final argument goes through the three following claims.

\begin{claim}\label{claim31} Let $(u_k)_{k\in\nn}$ be as in Proposition~$\ref{prop:asymp}$. Then
$$\limsup_{k\to \infty}\int_{\rn}(u_k)_+^{\crit}\, dx>0.$$\end{claim}

\noindent{\it Proof of Claim~$\ref{claim31}$:} We argue by
contradiction and assume that
\begin{equation}\label{hyp:contradic}
\lim_{k\to \infty}\int_{\rn}(u_k)_+^{\crit}\, dx=0.
\end{equation}
Estimating $\langle\Phi'(u_k), u_k\rangle$ and using inequality
\eqref{def:Ks} and \eqref{hyp:contradic}, we get as $k\to\infty$
\begin{eqnarray}
&&\Vert u_k\Vert^p=\Vert (u_k)_+\Vert_{L^{\crits}(\rn, |x|^{-s})}^{\crits}+o(1),\nonumber\\
&&\Vert (u_k)_+\Vert^{p}_{L^{\crits}(\rn, |x|^{-s})}\leq K(n,p,\mu,s)\Vert (u_k)_+\Vert_{L^{\crits}(\rn, |x|^{-s})}^{\crits}+o(1),\nonumber\\
&&\Vert (u_k)_+\Vert^{p}_{L^{\crits}(\rn, |x|^{-s})}\!
\left(1-K(n,p,\mu,s)\Vert (u_k)_+\Vert_{L^{\crits}(\rn,
|x|^{-s})}^{\crits-p}\right)\leq o(1). \label{ineq:contradic}
\end{eqnarray}
As in \eqref{balance:c} and \eqref{def:abc}, we have that
$$\int_{\rn}\frac{(u_k)_+^{\crits}}{|x|^s}\, dx= \frac{c\,p(n-s)}{p-s}+o(1)
$$
as $k\to \infty$. Plugging this inequality in \eqref{ineq:contradic}
and using the upper bound \eqref{bnd:c} on $c$, we get that
$$\lim_{k\to \infty}\int_{\rn}\frac{(u_k)_+^{\crits}}{|x|^s}\, dx=0.$$
A contradiction with \eqref{balance:c} and \eqref{hyp:contradic}
since $c>0$. This proves Claim~\ref{claim31}. \hfill$\Box$

\begin{claim}\label{claim32} Let $(u_k)_{k\in\nn}$ be a sequence as in Proposition~$\ref{prop:asymp}$.
Then there exists $\epsilon_1\in(0,\epsilon_0/2]$, with $\epsilon_0$
given in \eqref{4}, such that for all $\epsilon\in (0,\epsilon_1)$,
there exists a sequence $(r_k)_{k\in\nn}$ of $\rr_{>0}$ such that the
sequence $(\tilde{u}_k)_{k\in\nn}$ of $\dunp$, defined by
$$\tilde{u}_k(x):=r_k^{(n-p)/{p}}u_k(r_k x)\quad\mbox{ for }x\in\rn,$$
is again a Palais--Smale sequence of type given in Proposition~$\ref{prop:asymp}$ and verifies
\begin{equation}\label{hyp:epsilon}
\int_{B_1(0)}(\tilde{u}_k)_+^{\crit}\, dx=\epsilon
\end{equation}
for all $k\in\nn$.\end{claim}

\noindent{\it Proof of Claim~$\ref{claim32}$:} Let
$\lambda:=\limsup_{k\to \infty}\int_{\rn} (u_k)_+^{\crits} dx$. It
follows from Claim~\ref{claim31} that $\lambda>0$. Let
$\epsilon_1:=\min\{\epsilon_0/2,\lambda\}$, with $\epsilon_0>0$
given in \eqref{4}, see also Proposition~\ref{prop:ps}, and  fix
$\epsilon\in (0,\epsilon_1)$. Up to a subsequence, still denoted by
$(u_k)_{k\in\nn}$, for any $k\in\nn$ there exists $r_k>0$ such that
$$\int_{B_{r_k}(0)}(u_k)_+^{\crit}\, dx=\epsilon.$$
Due to scaling invariance, it is then straightforward to check that
$(\tilde{u}_k)_{k\in\nn}$ satisfies \eqref{hyp:epsilon} and the
properties  of Proposition~\ref{prop:asymp}. This proves
Claim~\ref{claim32}.\hfill$\Box$

\begin{claim}[Proof of Theorem \ref{th:main} when $\mu\ge0$]\label{claim33} Let $\tilde{u}_\infty\in\dunp$ be the weak limit
of $(\tilde{u}_k)_{k\in\nn}$ as $k\to \infty$ $($after a
subsequence$)$. Then $\tilde{u}_\infty\in C^1(\rn\setminus\{0\})$,
$\tilde{u}_\infty>0$ in $\rn\setminus\{0\}$ and $\tilde{u}_\infty$
is a weak solution of \eqref{eq:th}.
\end{claim}

\medskip\noindent{\it Proof of Claim~$\ref{claim33}$:}
We first assert that $(\tilde u_k)_k$ is bounded in  $D_1^p(\rn)$.
Indeed, since $p<\crit<\crits$ and $(\tilde u_k)_k$ is a
Palais--Smale sequence, there exist two positive constants $c_1$
and $c_2$ such that
$$\begin{aligned}c_1+c_2\|\tilde u_k\|&\ge\Phi(\tilde
u_k)-\frac1\crits\langle\Phi'(\tilde u_k),\tilde u_k\rangle\\
&=\left(\frac1p-\frac1{\crits}\right)\|\tilde u_k\|^p
+\left(\frac1\crit-\frac1\crits\right)\|(\tilde u_k)_+\|^{\crit}_{\crit}\\
&\ge\left(\frac1p-\frac1{\crits}\right)\|\tilde u_k\|^p
\end{aligned}$$
and the assertion follows at once by \eqref{comp:norm}, being $p>1$. Let $\tilde{u}_\infty\in\dunp$ be the weak limit
of $(\tilde{u}_k)_{k\in\nn}$ as $k\to \infty$, up to  a
subsequence. In case $\tilde{u}_\infty\equiv 0$, Proposition~\ref{prop:asymp} yields that either we have that $\lim_{k\to \infty}\int_{B_1(0)}(\tilde{u}_k)_+^{\crit}\, dx=0$ or we have that $\limsup_{k\to \infty}\int_{B_1(0)}(\tilde{u}_k)_+^{\crit}\, dx\ge\varepsilon_0$. Since $0<\varepsilon<\varepsilon_0/2$, this is a contradiction with \eqref{hyp:epsilon}.
Then $\tilde{u}_\infty\not\equiv 0$. It follows from Evans
\cite{evans} and Demengel--Hebey \cite{dh} (Lemmae~2 and~3) (see
also Saintier \cite{saintier1} Step~1.2 on p.303) that
$\tilde{u}_\infty$ is a nontrivial weak solution of
\begin{equation}\label{eq:tu}
-\Delta_p
\tilde{u}_\infty-\mu\frac{|\tilde{u}_\infty|^{p-2}\tilde{u}_\infty}{|x|^p}=(\tilde{u}_\infty)_+^{\crit-1}+\frac{(\tilde{u}_\infty)_+^{\crits-1}}{|x|^s}\quad\mbox{
in }\rn.
\end{equation}
We write \eqref{eq:tu} as
$-\Delta_p\tilde{u}_\infty=f(x,\tilde{u}_\infty)$, with an obvious
choice of $f$. Indeed, for all
$\omega\subset\subset\rn\setminus\{0\}$, there exists $C(\omega)>0$
such that $|f(x,u)|\leq C(\omega)(1+|u|^{\crit-1})$ for all
$x\in\omega$ and $u\in\rr$: it then follows from Theorem~2.1 of
Pucci--Servadei \cite{pucciservadei} (see also Druet
\cite[Lemmas~2.1 and~2.2]{druet}, Guedda--Veron
\cite[Proposition~1.1]{gv}) that $\tilde{u}_\infty\in
L^\infty_\text{\scriptsize{loc}}(\rn\setminus\{0\})$. Hence it
follows from Tolksdorf \cite[Theorem~1]{tolk} that
$\tilde{u}_\infty\in C^1(\rn\setminus\{0\})$.\par \noindent
Multiplying \eqref{eq:tu} by $(\tilde{u}_\infty)_-$ and integrating,
we get that $\Vert (\tilde{u}_\infty)_-\Vert=0$, and therefore
$(\tilde{u}_\infty)_-\equiv 0$ thanks to \eqref{comp:norm}. It then
follows that $\tilde{u}_\infty\in C^1(\rn\setminus\{0\})$ is a
non--negative nontrivial weak solution to \eqref{eq:tu}: thus
$\tilde{u}_\infty>0$ by the strong maximum principle of V\`azquez
\cite{vaz}. Therefore, $\tilde{u}_\infty\in\dunp\cap
C^1(\rn\setminus\{0\})$ is a {\it positive weak solution} of
\eqref{eq:th}. This proves Claim~\ref{claim33} and therefore
Theorem~\ref{th:main}.\hfill$\Box$

\medskip\noindent{\it Remark:} Consider the functional
$$\tilde{\Phi}(u):=\frac{1}{p}\Vert u\Vert^p-\frac{1}{\crit}\int_{\rn}|u|^{\crit}\, dx-\frac{1}{\crits}\int_{\rn}\frac{|u|^{\crits}}{|x|^s}\, dx$$
for $u\in \dunp$. Then the analysis above can be carried out for the
functional $\tilde{\Phi}$, with only minor modifications. The main
difference here is that the weak limit $\tilde{u}_\infty$ is not
necessarily positive.

\section{A non--existence result}\label{sec:poho}
In this section we require only that $\mu<\mu_1$ and prove the
following result:

\begin{thm}\label{th:nonex} Let $1<p<n$.
If $u\in \dunp$ is a weak solution to
\begin{equation}\label{eq:nonex}
-\Delta_p
u-\mu\frac{|u|^{p-2}u}{|x|^p}=\frac{|u|^{\crits-2}u}{|x|^s}+|u|^{q-2}u\quad\hbox{
in }\rn,
\end{equation}
when $s\in(0,p)$ and $1<q<\crit$, then $u\equiv 0$.
\end{thm}

\medskip\noindent{\it Remark 1}: note that, since $1<q<\crit$, we get that $u\in
L^q_\text{\scriptsize{loc}}(\rn)$ and the definition of the weak
solution makes sense.\par
\smallskip\noindent{\it Remark 2}: when $q>\crit$, the same conclusion holds if
$u\in L^\infty_\text{\scriptsize{loc}}(\rn\setminus\{0\})$ (see
Claims~\ref{claimp4} and~\ref{claimp5}).

\bigskip The proof of Theorem~\ref{th:nonex} uses a Pohozaev--type identity. It proceeds in five claims:

\begin{claim}\label{claimp1} Let $\eta$, $u\in C^\infty_c(\rn)$. Then

\begin{equation}\label{eq:claimp1}
\int_{\rn}|\nabla u|^{p-2}(\nabla u,\nabla(x,\nabla(\eta u)))\,
dx+\frac{n-p}{p}\int_{\rn}\eta|\nabla u|^p\, dx =B(u,\eta),
\end{equation}
where
$$\begin{aligned}
B(u,\eta)&= \int_{\rn} \Biggl(u|\nabla u|^{p-2}(\nabla u,\nabla\eta)+\nabla^2\eta(x,\nabla u)|\nabla u|^{p-2}u\Biggr.\\
&\qquad\qquad\qquad\Biggl.+|\nabla u|^{p-2}(\nabla
u,\nabla\eta)(x,\nabla u) +\frac{1}{p'}(x,\nabla\eta)|\nabla
u|^p\Biggr)dx,
\end{aligned}$$
and $p'=p/(p-1)$.\end{claim}

\noindent{\it Proof of Claim~$\ref{claimp1}$:} A similar identity
was proved by Guedda--Veron \cite{gv} on bounded domains of $\rn$.
Expanding $\nabla(x,\nabla(\eta u))$, we obtain that
\begin{eqnarray}
&&\int_{\rn}\!|\nabla u|^{p-2}(\nabla u,\nabla(x,\nabla(\eta u)))\,
dx\nonumber
=\int_{\rn}\!\eta|\nabla u|^p\, dx+\int_{\rn}\!\eta |\nabla u|^{p-2}x^i\partial_{ij}u\partial_j u\, dx\\
&&\,\,\,+ \int_{\rn} \Bigl(u|\nabla u|^{p-2}(\nabla
u,\nabla\eta)+\nabla^2\eta(x,\nabla u) |\nabla u|^{p-2}u
+|\nabla u|^{p-2}(\nabla u,\nabla\eta)(x,\nabla u)\Bigr.\label{eq:p1}\\
&&\,\,\,\Bigl.+(x,\nabla\eta)|\nabla u|^p\Bigr)dx,\nonumber
\end{eqnarray}
with Einstein's summation convention being used. Independently, we
have that
\begin{equation}\label{eq:p2}\begin{aligned}
\int_{\rn}\eta |\nabla u|^{p-2}x^i\partial_{ij}u\partial_j u\, dx&=\int_{\rn}\eta x^i\partial_i\left(\frac{|\nabla u|^p}{p}\right) dx\\
&= -\int_{\rn}\frac{\partial_i(\eta x^i)}{p}|\nabla u|^p
dx.\end{aligned}\end{equation} Plugging together \eqref{eq:p1} and
\eqref{eq:p2}, we get \eqref{eq:claimp1} and Claim~\ref{claimp1} is
proved.\hfill$\Box$

\begin{claim}\label{claimp2}
If $u\in D_1^p(\rn)\cap C^1(\rn\setminus\{0\})\cap
H_\text{\scriptsize{\rm 2,\,loc}}^1(\rn\setminus\{0\})$ and $\eta\in
C^\infty_c(\rn\setminus\{0\})$, then identity \eqref{eq:claimp1}
holds.
\end{claim}

\noindent{\it Proof of Claim~$\ref{claimp2}$:} By a density
argument, we get that there exists a sequence
$(\varphi_k)_{k\in\nn}\in C^\infty_c(\rn\setminus\{0\})$ such that
$\lim_{k\to \infty}\varphi_k=u$ in
$C^1_\text{\scriptsize{loc}}(\rn\setminus\{0\})\cap
H_\text{\scriptsize{\rm 2,\,loc}}^1(\rn\setminus\{0\})$. We then
apply Claim~\ref{claimp1} to $\eta,\varphi_k$ and let $k\to \infty$.
Claim~\ref{claimp2} is now proved. \hfill$\Box$

\begin{claim}\label{claimp3} Let  $f\in C^0((\rn\setminus\{0\})\times \rr)$ and let
$u\in D_1^p(\rn)\cap C^1(\rn\setminus\{0\}) \cap
H_\text{\scriptsize{\rm 2,\,loc}}^1(\rn\setminus\{0\})$ be a weak
solution of
\begin{equation}\label{eq:u:f}
-\Delta_p u=f(x,u)\quad\hbox{ in }\rn.
\end{equation}
Define $F(x,u):=\int_0^uf(x,v)\, dv$ and assume that $F\in
C^1((\rn\setminus\{0\})\times \rr)$. Moreover, along the solution
$u$,  assume that $u f(\cdot,u)$, $F(\cdot,u)$ and $x^i(\partial_i
F)(\cdot, u)\in L^1(\rn)$. Then
\begin{equation}\label{id:poho}
\int_{\rn}\left[\frac{n-p}{p}uf(x,u)-nF(x,u)-x^i(\partial_iF)(x,u)\right]
dx=0.
\end{equation}
\end{claim}

\noindent{\it Proof of Claim~$\ref{claimp3}$:} Fix $\eta\in
C^\infty_c(\rn\setminus\{0\})$. Using the notations of the proof of
Claim~\ref{claimp2} and~\eqref{eq:u:f}, we get that
\begin{equation}\label{eq:poho:1}\begin{aligned}
&\int_{\rn}|\nabla u|^{p-2}(\nabla u,\nabla(x,\nabla(\eta u)))\, dx\\
&\,\,\,=\lim_{k\to \infty}\int_{\rn}|\nabla u|^{p-2}(\nabla u,\nabla(x,\nabla(\eta \varphi_k)))\, dx\\
&\,\,\,=\lim_{k\to \infty} \int_{\rn}\! f(x,u) (x,\nabla(\eta\varphi_k))\, dx= \int_{\rn}\! f(x,u) (x,\nabla(\eta u))\, dx\\
&\,\,\,= \int_{\rn} u f(x,u) (x,\nabla \eta)\, dx
+\int_{\rn}\eta x^i\!\left[\partial_i(F(x,u))-(\partial_iF)(x,u)\right]dx \\
&\,\,\,= \int_{\rn} u f(x,u) (x,\nabla \eta)\, dx
-\int_{\rn}\partial_i(\eta x^i)F(x,u)\, dx-\int_{\rn}\eta
x^i(\partial_i F)(x,u)\, dx.
\end{aligned}\end{equation}
Independently, using \eqref{eq:u:f}, we have that
$$\begin{aligned}
\int_{\rn}|\nabla u|^{p-2}(\nabla u,\nabla(\eta u))\, dx&
=\lim_{k\to \infty}\int_{\rn}|\nabla u|^{p-2}(\nabla u,\nabla(\eta\varphi_k))\, dx\\
&=\lim_{k\to \infty}\int_{\rn} f(x,u)\eta\varphi_k\, dx=\int_{\rn}
f(x,u)\eta u\, dx
\end{aligned}$$
and therefore
\begin{equation}\label{eq:poho:2}
\int_{\rn}\eta|\nabla u|^p\, dx=\int_{\rn}\eta uf(x,u)\, dx
-\int_{\rn}u|\nabla u|^{p-2}(\nabla u,\nabla\eta)\, dx.
\end{equation}
Plugging \eqref{eq:poho:1} and \eqref{eq:poho:2} into
\eqref{eq:claimp1}, we get by H\"older's inequality that
\begin{equation}\label{eq:poho:finale}\begin{aligned}
&\left|\int_{\rn}\eta\left[\frac{n-p}{p}uf(x,u)-nF(x,u)-x^i(\partial_iF)(x,u)\right]\, dx\right|\\
&\qquad\qquad\leq \Vert\nabla
u\Vert^{p-1}_{L^p\left(\text{\scriptsize{\rm
Supp}}|\nabla\eta|\right)}\Vert u\Vert_{\crit}
\left(\frac np\Vert\nabla\eta\Vert_n+\Vert\, |x|\cdot|\nabla^2\eta|\,\Vert_{n}\right)\\
&\qquad\qquad+\Vert\, |x|\cdot|\nabla\eta|\,\Vert_\infty
\int_{\text{\scriptsize{\rm Supp}}|\nabla\eta|}
|uf(x,u)-F(x,u)|\, dx\\
&\qquad\qquad+\left(1+\frac{1}{p'}\right)\Vert\,
|x|\cdot|\nabla\eta| \,\Vert_\infty \Vert\nabla
u\Vert_{L^p\left(\text{\scriptsize{\rm Supp}}|\nabla\eta|\right)}.
\end{aligned}\end{equation}
We are left with choosing an appropriate cut--off function $\eta$.
Let $h\in C^\infty(\rr)$ be such that $h_{|\{t\leq 1\}}\equiv 0$,
$h_{|\{t\geq 2\}}\equiv 1$ and $0\leq h\leq 1$. Given $\epsilon>0$
small,  define $\eta_\epsilon$ as follows:
$\eta_\epsilon(x)=h(|x|/\epsilon)$ if $|x|\leq 3\epsilon$,
$\eta_\epsilon(x)=h(1/\epsilon |x|)$ if $|x|\geq (2\epsilon)^{-1}$
and $\eta_\epsilon(x)=1$ elsewhere. Clearly $\eta_\epsilon\in
C^\infty_c(\rn\setminus\{0\})$. Taking $\eta=\eta_\epsilon$ in
\eqref{eq:poho:finale} and letting $\epsilon\to 0$, we get
\eqref{id:poho} and Claim~\ref{claimp3} is proved. \hfill$\Box$

\begin{claim}\label{claimp4} If
$u\in D_1^p(\rn)\cap C^1(\rn\setminus\{0\})\cap
H_\text{\scriptsize{\rm 2,\,loc}}^1(\rn\setminus\{0\})$ is a weak
solution to \eqref{eq:nonex} when $q>1$ and $q\neq \crit$, then
$u\equiv 0$.
\end{claim}

\noindent{\it Proof of Claim~$\ref{claimp4}$:} In order to use
Claim~\ref{claimp3}, we need to prove that $u\in L^q(\rn)$. Indeed,
testing \eqref{eq:nonex} on $\eta_\epsilon u$, where
$\eta_\epsilon\in C_c^\infty(\rn\setminus\{0\})$ is as above (this
is a valid test--function, see the proof of \eqref{eq:poho:2}), we
get that
$$\int_{\rn}|\nabla u|^{p-2}(\nabla u,\nabla(\eta_\epsilon u))\, dx-\mu\int_{\rn}\frac{\eta_\epsilon|u|^p}{|x|^p}\, dx=\int_{\rn}\frac{\eta_\epsilon|u|^{\crits}}{|x|^s}\, dx+\int_{\rn}\eta_\epsilon|u|^q\, dx.$$
The Hardy inequality \eqref{ineq:hardy}, the Hardy--Sobolev
inequality \eqref{def:Ks} and H\"older's inequality yield  the
existence of $C>0$, independent of $\epsilon$, such that
$\int_{\rn}\eta_\epsilon |u|^q\, dx\leq C$ for all $\epsilon>0$.
Letting $\epsilon\to 0$, we get that $u\in L^q(\rn)$. Then we can
use Claim \ref{claimp3} and, applying \eqref{id:poho}, we have that
$$\left(\frac{1}{\crit}-\frac{1}{q}\right)\int_{\rn}|u|^q\, dx=0.$$
The fact that $q\neq\crit$ implies that $u\equiv 0$, and
Claim~\ref{claimp4} is proved. \hfill$\Box$

\begin{claim}\label{claimp5} Let $u\in D_1^p(\rn)$ be a weak solution of \eqref{eq:nonex},
with $q>1$.  Moreover, assume in addition that $u\in
L^\infty_{\text{\scriptsize{\rm loc}}}(\rn\setminus\{0\})$ in case
$q>\crit$. Then
$$u\in D_1^p(\rn)\cap C^1(\rn\setminus\{0\})\cap H_\text{\scriptsize{\rm 2,\,loc}}^1(\rn\setminus\{0\}).$$
\end{claim}

\noindent{\it Proof of Claim~$\ref{claimp5}$:} The argument relies
essentially on the works of Tolksdorf \cite{tolk}, Druet
\cite{druet} and Guedda--Veron \cite{gv}. We write \eqref{eq:nonex}
as $-\Delta_pu=f(x,u)$, with an obvious choice of $f$. Indeed, when
$1<q\leq \crit$, we get that for all
$\omega\subset\subset\rn\setminus\{0\}$, there exists $C(\omega)>0$
such that $-\Delta_pu=f(x,u)$, with $|f(x,u)|\leq
C(\omega)(1+|u|^{\crit-1})$ for all $x\in\omega$ and $u\in\rr$: it
then follows from Druet \cite[Lemmas~2.1 and~2.2]{druet},
Guedda--Veron \cite[Proposition~1.1]{gv} that $u\in
L^\infty_\text{\scriptsize{loc}}(\rn\setminus\{0\})$.

When also $u\in
L^\infty_\text{\scriptsize{loc}}(\rn\setminus\{0\})$, then $u$
satisfies $-\Delta_p u=f(x,u)$ weakly in $\rn$, with $f(\cdot,u)\in
L^\infty_\text{\scriptsize{loc}}(\rn\setminus\{0\})$. Hence it
follows from Tolksdorf \cite[Theorem~1 and
Pro\-po\-si\-tion~1]{tolk} that $u\in C^1(\rn\setminus\{0\})\cap
H_\text{\scriptsize{\rm 2,\,loc}}^{\inf\{2,p\}}(\rn\setminus\{0\})$.
This proves Claim~\ref{claimp5}. \hfill$\Box$

\medskip\noindent{\bf Proof of Theorem \ref{th:nonex}:}
The proof follows from the combination of Claims~\ref{claimp4}
and~\ref{claimp5}.

\section{Appendix 1: Extremals for Sobolev--type inequalities}\label{sec:app}
In this section we allow $\mu$ to be possible negative.

\begin{thm}\label{th:ext}
Let $p\in (1,n)$, $\mu<\mu_1$ and
 $s\in [0,p)$. If $s=0$, we assume that $\mu\geq 0$. Then the infimum $ {1}/{K(n,p,\mu,s)}$
in \eqref{def:Ks} is achieved. More precisely, if $(u_k)_{k\in\nn}$
is a minimizing sequence for ${1}/{K(n,p,\mu,s)}$ in $\dunp$ such
that $\int_{\rn}{|u_k|^{\crits}}|x|^{-s}dx=1$, then there exists a
sequence $(r_k)_{k\in\nn}$ in $\mathbb{R}_{>0}$ such that
$(r_k^{(n-p)/{p}}u_k(r_k\cdot))_{k\in\nn}$ is relatively compact in
$\dunp$ and converges to a minimizer for ${1}/{K(n,p,\mu,s)}$ up to
a subsequence. Moreover, the infimum is achieved by a non--negative
extremal.

Finally, if  $\mu\in[0,\mu_1)$ and if $s\in(0,p)$ when $\mu=0$, then any non--negative minimizer of \eqref{def:Ks}
in   $\dunp\setminus\{0\}$  is positive, radially symmetric,
radially decreasing with respect to $0$ and approaches zero as $|x|\to\infty$.
\end{thm}

\medskip\noindent{\bf Remark 1:} The assumption that $\mu\ge0$ in case $s=0$ is not technical.
Indeed, as shown in Claim~\ref{sec:neg:mu}.1, it is not difficult to
prove that $K(n,p,\mu,0)=K(n,p,0,0)$ when $\mu <0$: then, since
there are extremals for $K(n,p,0,0)$, there is no extremal for
$K(n,p,\mu,0)$. We refer to Lions \cite{lions} for further
considerations on this phenomenon.

\medskip\noindent{\bf Remark 2:} When $p=2$, the statement of Theorem \ref{th:ext} is essentially contained in Catrina-Wang \cite{cw}. In particular, the assumption that $\mu\geq 0$ in the last assertion of the theorem is not technical: indeed, it follows from Catrina-Wang \cite{cw} that when $p=2$, for any $\mu<0$, there exists $s_\mu>0$ such that for all $s\in (0,s_\mu)$, then no minimizer of \eqref{def:Ks} is radially symmetrical.

\bigskip The proof of Theorem \ref{th:ext} relies essentially on Lions's proof of the existence
of extremals for the classical Sobolev inequalities \cite{lions}. We
mainly follow the proof given in the book of Struwe \cite{struwe}.
Note that when $s=\mu=0$, the extremals exist (see Rodemich
\cite{rodemich}, Aubin \cite{aubin}, Talenti \cite{talenti}, see
also Lions \cite{lions}).

\medskip  Let $(\tilde{u}_k)_{k\in\nn}\subset\dunp\setminus\{0\}$ be a minimizing sequence
for $1/K(n,p,\mu,s)$ in \eqref{def:Ks}. Up to multiplying by a
positive constant, we assume that
\begin{equation*}
\int_{\rn}\frac{|\tilde{u}_k|^{\crits}}{|x|^s}\, dx=1
\qquad\mbox{and}\qquad \lim_{k\to \infty}\int_{\rn}\left(|\nabla
\tilde{u}_k|^p-\mu\frac{|\tilde{u}_k|^p}{|x|^p}\right)\,
dx=\frac{1}{K(n,p,\mu,s)}.
\end{equation*}
Since $\int_{\rn}{|\tilde{u}_k|^{\crits}}|x|^{-s} dx=1$ for all
$k\in\nn$, there exists $r_k>0$ such that
$$\int_{B_{r_k(0)}}\frac{|\tilde{u}_k|^{\crits}}{|x|^s}\,dx=\frac{1}{2}$$
for all $k\in\nn$. We define the rescaled sequence
$$u_k(x):=r_k^{(n-p)/{p}}\tilde{u}_k(r_k x)$$
for all $k\in\nn$ and $x\in\rn$. Clearly $u_k\in\dunp$ for all
$k\in\nn$ and $(u_k)_{k\in\nn}$ is a minimizing sequence for
$1/K(n,p,\mu,s)$, that is
\begin{equation}\label{ext:1}
\int_{\rn}\frac{|u_k|^{\crits}}{|x|^s}\, dx=1 \quad\mbox{and}\quad
\lim_{k\to \infty}\int_{\rn}\left(|\nabla
u_k|^p-\mu\frac{|u_k|^p}{|x|^p}\right)dx=\frac{1}{K(n,p,\mu,s)}.
\end{equation}
Moreover, we have that
\begin{equation}\label{eq:contradic:ext}
\int_{B_1(0)}\frac{|{u}_k|^{\crits}}{|x|^s}\, dx=\frac{1}{2}
\end{equation}
for all $k\in\nn$. In addition, $\Vert
u_k\Vert^p=K(n,p,\mu,s)^{-1}+o(1)$ as $k\to \infty$, and then, using
\eqref{comp:norm}, the $(\Vert\nabla u_k\Vert_p)_{k\in\nn}$ is
bounded. Therefore, without loss of generality, we assume that there
exists $u\in\dunp$ such that
\begin{eqnarray}
&&u_k\rightharpoonup u\quad\hbox{ weakly in }\dunp\hbox{ as }k\to \infty,\nonumber\\
&& \lim_{k\to \infty}u_k(x)=u(x)\quad\hbox{ for a.a.
}x\in\rn.\nonumber
\end{eqnarray}
We define the measures
\begin{equation}\label{vl}
\nu_k:=\frac{|u_k|^{\crits}}{|x|^s}\, dx\quad\hbox{ and }\quad
\lambda_k:=\left(|\nabla u_k|^p-\mu
\frac{|u_k|^{p}}{|x|^p}\right)dx.
\end{equation}
Hence \eqref{ext:1} simply reduces to
\begin{equation}\label{ext:1'}\int_{\rn}d\nu_k=1\quad\mbox{and}\quad
\lim_{k\to
\infty}\int_{\rn}d\lambda_k=\frac{1}{K(n,p,\mu,s)}.\end{equation}
Clearly, $\nu_k\geq 0$ by \eqref{ext:1}. Moreover, in the sense of
measures, we get that $|\lambda_k|\leq \left(|\nabla u_k|^p+|\mu|
{|u_k|^{p}}|x|^{-p}\right)dx$ is a bounded measure with respect to
$k\in\nn$. Up to a subsequence, there exist two measures $\nu$  and
$\lambda$ such that
$$\nu_k\rightharpoonup \nu\quad\hbox{ and }\quad\lambda_k\rightharpoonup\lambda\quad\hbox{weakly in the sense of
measures as }k\to \infty.$$

\smallskip\noindent  We now apply Lions's first concentration--compactness Lemma \cite{lions} to the
sequence of measures $(\nu_k)_{k\in\nn}$. Indeed, up to a
subsequence, three situations can occur (cf. \cite[Lemma~1, page
39]{struwe}):

\medskip(a) ({\it Compactness}) There exists a sequence $(x_k)_{k\in\nn}$ in $\rn$ such that for any $\epsilon>0$
there exists $R_\epsilon>0$ for which
$$\int_{B_{R_\varepsilon}(x_k)}d\nu_k\geq 1-\epsilon\quad\hbox{ for all }k\in\nn\hbox{ large}.$$

\smallskip(b) ({\it Vanishing}) For all $R>0$  there holds
$$\lim_{k\to \infty}\left(\sup_{x\in\rn}\int_{B_R(x)}d\nu_k\right)=0.$$

\smallskip(c) ({\it Dichotomy}) There exists $\alpha\in (0,1)$ such that for any $\epsilon>0$
there exists $R_\epsilon>0$ and a sequence
$(x_k^\epsilon)_{k\in\nn}\in\rn$, with the following property: given
$R'>R_\epsilon$, there are non--negative measures $\nu_k^1$ and
$\nu_k^2$ such that
$$\begin{gathered}
0\leq \nu_k^1+\nu_k^2\leq \nu_k,\quad\hbox{Supp}(\nu_k^1)\subset
B_{R_\epsilon}(x_k^\epsilon),\quad
\hbox{Supp}(\nu_k^2)\subset \rn\setminus B_{R'}(x_k^\epsilon),\\
\nu_k^1=\nu_k\bigl.\bigr|_{B_{R_\epsilon}(x_k^\epsilon)},\qquad
\nu_k^2=\nu_k\bigl.\bigr|_{\rn\setminus B_{R'}(x_k^\epsilon)},\\
\limsup_{k\to
\infty}\left(\left|\alpha-\int_{\rn}d\nu_k^1\right|+\left|(1-\alpha)-\int_{\rn}d\nu_k^2\right|\right)\leq\epsilon.
\end{gathered}$$

\begin{claim}\label{claim:conc} Compactness $($point {\rm(a)}$)$ holds. In particular, we have that $\int_{\rn}d\nu=1$.
\end{claim}

\begin{proof} It follows from \eqref{eq:contradic:ext} that {\it Vanishing}, point (b), does not hold.
We argue by contradiction and assume that {\it Dichotomy} holds,
that is there exists $\alpha\in (0,1)$ such that (c) above holds.
Taking $\epsilon=(k+1)^{-1}$, we can assume that, up to a
subsequence, there exist sequences $(R_k)_{k\in\nn}$ in $\rr_{>0}$,
$(x_k)_{k\in\nn}$ in $\rn$ and two sequences of non--negative
measures, $(\nu_k^1)_{k\in\nn}$ and $(\nu_k^2)_{k\in\nn}$, such that
\begin{eqnarray}
&&0\leq \nu_k^1+\nu_k^2\leq \nu_k,\qquad \lim_{k\to \infty}R_k=\infty,\nonumber\\
&&\hbox{Supp}(\nu_k^1)\subset B_{R_k}(x_k),\qquad
\hbox{Supp}(\nu_k^2)\subset
\rn\setminus B_{2R_k}(x_k),\nonumber\\
&&\label{lim:alpha}\begin{gathered}
\nu_k^1=\nu_k\bigl.\bigr|_{B_{R_k}(x_k)},\qquad
\nu_k^2=\nu_k\bigl.\bigr|_{\rn\setminus B_{2R_k}(x_k)},\\
\lim_{k\to \infty}\int_{\rn}d\nu_k^1=\alpha\quad\hbox{ and }\quad
\lim_{k\to \infty}\int_{\rn}d\nu_k^2=1-\alpha. \end{gathered}
\end{eqnarray}
In particular, by \eqref{ext:1'}$_1$ and \eqref{lim:alpha}, we have
\begin{equation}\label{est:ring}
\lim_{k\to \infty}\int_{D_k}d\nu_k=0,\qquad
D_k:=B_{2R_k}(x_k)\setminus B_{R_k}(x_k).
\end{equation}
\medskip\noindent{\it Step \ref{claim:conc}.1:} We claim that
\begin{equation}\label{est:ring:2}
\lim_{k\to \infty}\int_{D_k}\frac{|u_k|^p}{|x|^p}\, dx=0.
\end{equation}
Indeed, by  H\"older's inequality, we get that
$$\begin{aligned}
\int_{D_k}\frac{|u_k|^p}{|x|^p}\,
dx&=\int_{D_k}\frac{1}{|x|^{p-{ps}/{\crits}}}
\left(\frac{|u_k|}{|x|^{{s}/{\crits}}}\right)^p\, dx\\
\!\!&\leq
\left(\int_{D_k}\!\!\left(\frac{1}{|x|^{p-{ps}/{\crits}}}\!\right)^{(n-s)/(p-s)}\!\!\!\!
dx\right)^{1-\frac{p}{\crits}}
\!\!\left(\int_{D_k}\!\! \frac{|u_k|^{\crits}}{|x|^s}\, dx\right)^{{p}/{\crits}}\\
\!\!&\leq  C \left(\int_{D_k}d\nu_k\right)^{{p}/{\crits}}.
 \end{aligned}$$
Therefore, \eqref{est:ring} yields \eqref{est:ring:2}, and the claim
is proved.

\medskip\noindent{\it Step \ref{claim:conc}.2:} Let $\varphi\in C_c^\infty(\rn)$ such that $0\leq\varphi\leq 1$, $\varphi_{|B_1(0)}\equiv 1$ and $\varphi_{|B_2(0)^c}\equiv 0$.
We define $\varphi_k(x):=\varphi(R_k^{-1}(x-x_k))$ for all $x\in\rn$
and all $k\in\nn$. By \eqref{lim:alpha}, \eqref{est:ring},
\eqref{def:Ks} and the fact that $p<\crits$, we get that
\begin{equation}\label{est:sobo}\begin{aligned}
1&= \left(\int_{\rn}\varphi_k^{\crits}d\nu_k^1+\int_{\rn}(1-\varphi_k)^{\crits}d\nu_k^2\right)^{{p}/{\crits}}+o(1)\\
&\leq
\left(\int_{\rn}\varphi_k^{\crits}d\nu_k^1\right)^{{p}/{\crits}}
+\left(\int_{\rn}(1-\varphi_k)^{\crits}d\nu_k^2\right)^{{p}/{\crits}}+o(1)\\
&\leq
\left(\int_{\rn}\varphi_k^{\crits}d\nu_k\right)^{{p}/{\crits}}
+\left(\int_{\rn}(1-\varphi_k)^{\crits}d\nu_k\right)^{{p}/{\crits}}+o(1)\\
&\leq K(n,p,\mu,s)\int_{\rn}\left(|\nabla(\varphi_k u_k)|^p-\mu\frac{|\varphi_k u_k|^p}{|x|^p}\right) dx\\
&\qquad+ K(n,p,\mu,s)\int_{\rn}\left(|\nabla((1-\varphi_k)
u_k)|^p-\mu\frac{|(1-\varphi_k) u_k|^p}{|x|^p}\right) dx+o(1).
\end{aligned}\end{equation}
\medskip\noindent{\it Step \ref{claim:conc}.3:} As shown in \eqref{io}, we shall prove that
\begin{equation}\label{est:phi}
\int_{\rn}|\nabla(\varphi_k u_k)|^p\,
dx=\int_{\rn}|\varphi_k|^p|\nabla  u_k|^p\, dx+o(1)
\end{equation}
as $k\to\infty$. Indeed,
$$\left||\nabla(\varphi_k u_k)|^p-|\varphi_k|^p|\nabla u_k|^p\right|\leq
C_p(|\varphi_k \nabla u_k|^{p-1}|u_k\nabla\varphi_k|+
|u_k\nabla\varphi_k|^p)$$ for all $k\in\nn$, which, integrated over
$\rn$, gives
$$
\left|\int_{\rn}\!|\nabla(\varphi_k u_k)|^p
dx-\int_{\rn}\!|\varphi_k|^p|\nabla  u_k|^p dx\right| \leq
C_p\!\int_{\rn}\!(|\varphi_k \nabla u_k|^{p-1}|u_k\nabla\varphi_k|+
|u_k\nabla\varphi_k|^p) dx.$$ By H\"older's inequality, and since
$\hbox{Supp\,}(\nabla\varphi_k)\subset D_k$, we get that
\begin{eqnarray*}
&&\left|\int_{\rn}|\nabla(\varphi_k u_k)|^p dx-\int_{\rn}|\varphi_k|^p|\nabla  u_k|^p dx\right|\\
&&\qquad\leq C_p\Vert\nabla u_k\Vert_p^{p-1}\left(\int_{\rn}|u_k\nabla\varphi_k|^p\, dx\right)^{1/p}+C_p \int_{\rn}|u_k\nabla\varphi_k|^p\, dx\\
&&\qquad\leq C\left[ \left(\int_{D_k}\!\!\frac{|u_k|^p}{|x|^p}\,
dx\right)^{1/p}+\int_{D_k}\frac{|u_k|^p}{|x|^p}\, dx\right]
\end{eqnarray*}
for all $k\in\nn$. Therefore, \eqref{est:ring:2} yields
\eqref{est:phi}. This ends Step~\ref{claim:conc}.3.

\medskip\noindent{\it Step \ref{claim:conc}.4:} Similarly to \eqref{est:phi} it results
\begin{equation}\label{est:1-phi}
\int_{\rn}|\nabla\{(1-\varphi_k) u_k\}|^p\,
dx=\int_{\rn}|1-\varphi_k|^p|\nabla  u_k|^p\, dx+o(1)
\end{equation}
as $k\to\infty$. Plugging \eqref{est:phi} and \eqref{est:1-phi} into
\eqref{est:sobo}, we obtain
\begin{equation}\label{ineq:3}\begin{aligned}
1&\leq
\left(\int_{\rn}\varphi_k^{\crits}d\nu_k^1\right)^{{p}/{\crits}}
+\left(\int_{\rn}(1-\varphi_k)^{\crits}d\nu_k^2\right)^{{p}/{\crits}}+o(1)\\
&\leq  K(n,p,\mu,s)\int_{\rn}\left[\varphi_k^p+(1-\varphi_k)^p\right] d\lambda_k+o(1)\\
&= 1
+K(n,p,\mu,s)\int_{\rn}\left[\varphi_k^p+(1-\varphi_k)^p-1\right]
d\lambda_k+o(1)
 \end{aligned}\end{equation}
by \eqref{ext:1'}$_2$. We now deal with the second term of the right
hand side above. Since
$\hbox{Supp\,}(1-\varphi_k^p-(1-\varphi_k)^p)\subset D_k$ and $0\leq
\varphi_k^p+(1-\varphi_k)^p\leq 1$, we get that
\begin{equation}\label{err:phi}\begin{aligned}
\int_{\rn}\left[\varphi_k^p+(1-\varphi_k)^p-1\right] d\lambda_k&=
-\int_{\rn}\left[1-\varphi_k^p-(1-\varphi_k)^p\right]|\nabla u_k|^p dx \\
&\,\,\,\,\,\,-\mu\int_{D_k}\!\!\!\!\!\!\left[\varphi_k^p+(1-\varphi_k)^p-1\right]\frac{|u_k|^p}{|x|^p}  dx \\
&\leq  2|\mu|\int_{D_k}\frac{|u_k|^p}{|x|^p}\, dx.
 \end{aligned}\end{equation}
Letting $k\to \infty$ in   \eqref{ineq:3} and using \eqref{err:phi},
\eqref{est:ring:2} and \eqref{lim:alpha}, we get that
$1=\alpha^{{p}/{\crits}}+(1-\alpha)^{{p}/{\crits}}$. This is
impossible when $\alpha\in (0,1)$, since $p<\crits$ being
$s\in[0,p)$. This contradiction proves
Claim~\ref{claim:conc}.\end{proof}

\begin{claim}\label{claim:compconc} There exist $J\subset \nn$ at most countable, a subset $I\subset  J$ and
a family $\{x_i\}_{i\in J}$ in $\rn$ such that
\begin{equation}\label{dec:nu}
\nu=\frac{|u|^{\crits}}{|x|^s}\, dx+\sum_{i\in I}\nu^i\delta_{x_i},
\end{equation}
where $\nu^i=\nu(\{x^i\})>0$ for all $i\in I$. In particular,
$\{x_i\,/\,i\in I\}\subset \{0\}$ when $s>0$. Moreover, there exists
a bounded non--negative measure $\lambda_0\geq 0$ with no atoms
$($that is $\lambda_0(\{x\})=0$ for all $x\in\rn)$ such that
\begin{equation}\label{dec:lambda}
\lambda=\lambda_0+\left(|\nabla u|^p-\mu\frac{|u|^p}{|x|^p}\right)\,
dx+\sum_{i\in J}\lambda^i\delta_{x_i},
\end{equation}
where $\lambda^i=\lambda(\{x_i\})>0$ for all $i\in J$. In addition,
\begin{equation}\label{comp:atoms}
(\nu^i)^{{p}/{\crits}}\leq K(n,p,\mu,s)\lambda^i\quad\hbox{ for all
}i\in I.
\end{equation}
\end{claim}
\begin{proof} This proof is essentially an adaptation of Lions's second concentration--compactness Lemma \cite{lions}.
When $s=0$, \eqref{dec:nu} is a consequence of Lions's result. When
$s>0$, since $(u_k)_{k\in\nn}$ goes to $u$ strongly in
$L^q_{\text{\scriptsize{\rm loc}}}(\rn)$ for $q<\crit$, we get that
$\nu={|u|^{\crits}}|x|^{-s}\, dx+\nu(\{0\})\delta_0$. This proves
\eqref{dec:nu} in the case $s\geq 0$.\par

\medskip\noindent We are left with proving \eqref{dec:lambda}. As above, we get that there exists $L\geq 0$ such that
\begin{equation}\label{dec:hardy}
\frac{|u_k|^p}{|x|^p}\, dx\rightharpoonup \frac{|u|^p}{|x|^p}\,
dx+L\,\delta_0
\end{equation}
in the sense of measures as $k\to \infty$. Up to extraction, we let
$\lambda'$ be the weak limit of $(|\nabla u_k|^p\, dx)$ as $k\to
\infty$ in the sense of measures. Since $u_k\rightharpoonup u$
weakly in $\dunp$ as $k\to \infty$, we get that $\lambda'\geq
|\nabla u|^p\, dx$. Therefore, we decompose $\lambda'$ as follows:
\begin{equation}\label{dec:lambdaprime}
\lambda'=\lambda_0+|\nabla u|^p\, dx+\sum_{j\in
K}\lambda'(\{z_j\})\delta_{z_j},
\end{equation}
where $\lambda_0\geq 0$ and the $z_j$'s, $j\in K$ countable, are the
atoms of $\lambda'$. Combining \eqref{dec:hardy} and
\eqref{dec:lambdaprime}, we have that
\begin{equation}\label{dec:lambda:1}
\lambda=\lambda_0+\left(|\nabla u|^p-\mu\frac{|u|^p}{|x|^p}\right)\,
dx -L\mu\,\delta_0 +\sum_{j\in K}\lambda'(\{z_j\})\delta_{z_j}.
\end{equation}
\medskip\noindent We claim that
\begin{equation}\label{comp:atom}
[\nu(\{x\})]^{{p}/{\crits}}\leq
K(n,p,\mu,s)\lambda(\{x\})\quad\hbox{ for all }x\in\rn.
\end{equation}
\noindent Indeed, take $\varphi\in C_c^\infty(\rn)$ such that
$\varphi_{|B_1(0)}\equiv 1$ and $\varphi_{|\rn\setminus
B_2(0)}\equiv 0$. Given $x_0\in\rn$ and $\epsilon>0$, we define
$\varphi_{\epsilon}(x)=\varphi(\epsilon^{-1}(x-x_0))$ for $x\in\rn$.
It follows from the Sobolev inequality \eqref{def:Ks} that
$$\left(\int_{\rn}\frac{|\varphi_\epsilon u_k|^{\crits}}{|x|^s}\, dx\right)^{{p}/{\crits}}\leq
K(n,p,\mu,s)\int_{\rn}\left(|\nabla (\varphi_\epsilon
u_k)|^p-\mu\frac{|\varphi_\epsilon u_k|^p}{|x|^p}\right)\, dx.$$ As
in the proof of \eqref{est:phi}, we have
\begin{eqnarray}\qquad
\left(\int_{\rn}|\varphi_\epsilon|^{\crits}\,
d\nu_k\right)^{{p}/{\crits}}\leq
K(n,p,\mu,s)\int_{\rn}|\varphi_\epsilon|^p\, d\lambda_k+
C\theta_k+C(\theta_k)^{1/p}
\end{eqnarray}
for all $k\in\nn$ and all $\epsilon>0$,  where
$$\theta_k:=\int_{B_{2\epsilon}(x_0)\setminus B_{\epsilon}(x_0)}\frac{|u_k|^p}{|x|^p}\, dx.$$
Letting $k\to \infty$ and then $\epsilon\to 0$, we get that
$$[\nu(\{x_0\})]^{{p}/{\crits}}\leq K(n,p,\mu,s)\lambda(\{x_0\})$$
and the claim is proved.

\medskip Combining \eqref{dec:lambda:1} with \eqref{comp:atom} and considering separately the cases $0\in\{x_i/\, i\in J\}$ or not, we get \eqref{dec:lambda}. This proves Claim \ref{claim:compconc}.\end{proof}

\begin{claim}\label{claim:dicho} We assert that
$$\hbox{either }\left\{\nu=\frac{|u|^{\crits}}{|x|^s}\, dx\hbox{ and }\lambda=
\left(|\nabla u|^p-\mu\frac{|u|^p}{|x|^p}\right) dx\right\}$$
$$\hbox{or there exists $x_0\in\rn$ such that } \left\{\nu=\delta_{x_0}\hbox{ and }
\lambda=\frac{\delta_{x_0}}{K(n,p,\mu,s)}\right\}.$$
\end{claim}
\begin{proof}
Integrating \eqref{dec:nu} and \eqref{dec:lambda}, using
\eqref{def:Ks}, \eqref{comp:atoms} and the fact that
$\int_{\rn}d\nu=1$ (see Claim~\ref{claim:conc}) and inequality
\eqref{def:K}, we have
\begin{equation}\begin{aligned}\label{ineq:lambda}
1&=\left(\int_{\rn}d\nu\right)^{p/\crits}=\left(\int_{\rn}\frac{|u|^{\crits}}{|x|^s}\, dx+\sum_{i\in I}\nu^i\right)^{p/\crits}\\
&\leq \left(\int_{\rn}\frac{|u|^{\crits}}{|x|^s}dx\right)^{p/\crits}+\sum_{i\in I}(\nu^i)^{p/\crits}\\
&\leq  K(n,p,\mu,s)\left(\int_{\rn}\left(|\nabla u|^p-\mu\frac{|u|^p}{|x|^p}\right) dx+\sum_{i\in I}\lambda^i\right)\\
&\leq  K(n,p,\mu,s)\int_{\rn}d\lambda.
\end{aligned}\end{equation}
We are then left with estimating $\int_{\rn}d\lambda$ from above.
Let $\psi\in C^\infty(\rn)$ such that $0\leq\psi\leq 1$,
$\psi_{|B_1(0)}\equiv 0$ and $\psi_{|\rn\setminus B_2(0)}\equiv 1$.
Given $R>0$, we let
 $\psi_R(x)=\psi(R^{-1}x)$ for $x\in\rn$. In particular, $1-\psi_R^p\in C^0_c(\rn)$. Hence, since $\mu<\mu_1$,
by \eqref{ext:1'}$_2$ and \eqref{ineq:hardy} we find that
\begin{equation}\begin{aligned}\label{lim:psiR}
\int_{\rn} (1-\psi_R^p)\, d\lambda_k&=\int_{\rn}d\lambda_k-\int_{\rn}\left(\psi_R^p|\nabla u_k|^p-\mu\frac{|\psi_R u_k|^p}{|x|^p}\right)\, dx\\
&= \int_{\rn}d\lambda_k-\int_{\rn}\left(|\nabla (\psi_R u_k)|^p-\mu\frac{|\psi_R u_k|^p}{|x|^p}\right)\, dx\\
&\qquad\qquad +\int_{\rn}\left(|\nabla (\psi_R u_k)|^p-\psi_R^p|\nabla u_k|^p\right)\, dx\\
&\leq \frac{1}{K(n,p,\mu,s)}+\int_{\rn}\left(|\nabla (\psi_R
u_k)|^p-\psi_R^p|\nabla u_k|^p\right)\, dx+o(1).
\end{aligned}\end{equation}
Mimicking what was worked out in \eqref{est:phi}, we obtain
\begin{equation*}
\left|\int_{\rn}\left(|\nabla (\psi_R u_k)|^p-\psi_R^p|\nabla
u_k|^p\right)\, dx\right|\leq  C\theta_{k}(R)+C\theta_k(R)^{p},
\end{equation*}
where
$$\theta_k(R):=\int_{B_{2R}(0)\setminus B_{R}(0)}\frac{|u_k|^p}{|x|^p}\, dx.$$
Therefore, letting $k\to \infty$ in \eqref{lim:psiR}, and then $R\to
\infty$, we get that
$$\int_{\rn}  d\lambda\leq  \frac{1}{K(n,p,\mu,s)}.$$
Plugging this latest inequality in \eqref{ineq:lambda}, we get that
$\int_{\rn}  d\lambda =  K(n,p,\mu,s)^{-1}$. Therefore, there is
equality in \eqref{ineq:lambda}. By convexity, this means that one
and only one term in \eqref{dec:nu} is nonzero and that there is
equality in all the inequalities used to prove \eqref{ineq:lambda}.
The conclusion of the claim then follows. This proves
Claim~\ref{claim:dicho}.
\end{proof}

\begin{claim}\label{claim:nonconc} We assert that $\nu={|u|^{\crits}}|x|^{-s}\, dx$ and
$\lambda=(|\nabla u|^p-\mu{|u|^p}|x|^{-p})\, dx$.
\end{claim}
\begin{proof} We argue by contradiction. If Claim~\ref{claim:nonconc} does not hold, it follows from Claim \ref{claim:dicho}, that there exists $x_0\in\rn$ such that
$\nu=\delta_{x_0}$ and $\lambda=\delta_{x_0}/{K(n,p,\mu,s)}$: in
particular, $u\equiv 0$. If $x_0=0$, then
 $\int_{B_{1/2}(0)}d\nu=1$, which contradicts the initial hypotheses \eqref{eq:contradic:ext}
 and proves Claim~\ref{claim:nonconc} when $x_0=0$. We are then left with proving that $x_0=0$. We argue by contradiction and assume that $x_0\neq 0$. We distinguish two cases:

\smallskip\noindent{\it Case 1:}  $s>0$. Then, since $u\equiv 0$, we get that
$\lim_{k\to \infty}u_k=0$ in $L^{\crits}_{\text{\scriptsize{\rm
loc}}}(\rn)$, and then $\lim_{k\to
\infty}\int_{B_\delta(x_0)}{|u_k|^{\crits}}|x|^{-s}\, dx=0$ for
$\delta>0$ small enough: a contradiction with the fact that
$\nu=\delta_{x_0}$. This ends Case~1.\par

\medskip\noindent{\it Case 2:}  $s=0$. Let $\delta>0$ and
$\varphi\in C_c^\infty(\rn)$ such that
$0\leq\varphi\leq\varphi(x_0)=1$ and $\varphi_{\rn\setminus
B_\delta(x_0)}\equiv 0$. Since $\lim_{k\to \infty}u_k=0$ in
$L^p_{\text{\scriptsize{\rm loc}}}(\rn)$, it follows from the Hardy
inequality \eqref{ineq:hardy} and computations similar to the ones
leading to \eqref{est:1-phi} that there exists $C>0$ such that
$$\begin{aligned}
\int_{\rn}\frac{|(1-\varphi)u_k|^p}{|x|^p}\, dx&
\leq C\int_{\rn}\left(|\nabla\{ (1-\varphi)u_k\}|^p-\mu\frac{|(1-\varphi)u_k|^p}{|x|^p}\right)dx\\
&=  C\int_{\rn}(1-\varphi)^pd\lambda_k+o(1)\\
&=  C\left(\frac{1}{K(n,p,\mu,s)}-\int_{\rn}[1-(1-\varphi)^p]d\lambda_k\right)+o(1)\\
&=
\frac{C}{K(n,p,\mu,s)}\left\{1-[1-(1-\varphi(x_0))^p]\right\}+o(1)=
o(1)
\end{aligned}$$
as $k\to \infty$, since clearly \eqref{est:1-phi} holds when
$\varphi$ replaces $\varphi_k$, and $\lambda_k\rightharpoonup
\lambda=\delta_{x_0}/{K(n,p,\mu,s)}$. In particular, for all
$\delta>0$, we get that
$$\lim_{k\to \infty}\int_{\rn\setminus B_\delta(x_0)}\frac{|u_k|^p}{|x|^p}\, dx=0.$$
Moreover, since $x_0\neq 0$ and $u_k\to 0$ strongly in
$L^p_\text{\scriptsize{loc}}(\rn)$, we have
$$\lim_{k\to \infty}\int_{\rn}\frac{|u_k|^p}{|x|^p}\, dx=0,$$
which implies by \eqref{ext:1'}, since $s=0$, that
$$\frac{\int_{\rn}|\nabla u_k|^p\, dx}{\left(\int_{\rn}|u_k|^{\crit}\, dx\right)^{{p}/{\crit}}}=\frac{1}{K(n,p,\mu,0)}+o(1)$$
as $k\to \infty$. It then follows from \eqref{def:K} that
\begin{equation}\label{ineq:K}
\frac{1}{K(n,p,0,0)}\leq \frac{1}{K(n,p,\mu,0)}.
\end{equation}
Let $u\in\dunp\setminus\{0\}$ be an extremal for $K(n,p,0,0)$ (this
exists, see Rodemich \cite{rodemich}, Talenti \cite{talenti}, Aubin
\cite{aubin} and also Lions \cite{lions}). Estimating the functional
of $K(n,p,\mu,0)$ at $u$ and using that $\mu>0$, we get that
$$\frac{1}{K(n,p,0,0)}> \frac{1}{K(n,p,\mu,0)}.$$
A contradiction with inequality \eqref{ineq:K}. This rules out the
case $x_0\neq 0$, and Case~2 is finished. This also ends the proof
of Claim \ref{claim:nonconc}.\end{proof}

\bigskip\noindent{\bf Proof of Theorem \ref{th:ext}.}
Since $\nu={|u|^{\crits}}|x|^{-s}dx$ and $\lambda=(|\nabla u|^p-\mu
\,{|u|^p}|x|^{-p}) dx$, we get that $\lim_{k\to \infty}u_k=u$ in
$L^{\crits}(\rn,|x|^{-s})\cap L^p(\rn,|x|^{-p})$. Consequently, we
get that $\Vert \nabla u_k\Vert_p\to\Vert \nabla u\Vert_p$ as $k\to
\infty$ and by Clarkson's uniform convexity, we find that
$\lim_{k\to \infty}u_k=u$ in $\dunp$. Hence $u$ is an extremal for
\eqref{def:Ks}. In addition, $|u|$ is in $\dunp$ and
$|\nabla|u||=|\nabla u|$ a.e on $\rn$: therefore, $|u|$ is also an
extremal, and then there exist non--negative extremals.
The first part of Theorem~\ref{th:ext} is proved.

Assume now that $\mu\in[0,\mu_1)$ and $s\in(0,p)$ when $\mu=0$. Let $u\ge0$ be a
minimizer of \eqref{def:Ks} in $\dunp\setminus\{0\}$, which exists
from the first part of Theorem~\ref{th:ext} already proved.
Following Talenti \cite{talenti}, see also \cite[Section~3.2]{ll},  we define the Schwarz symmetrization of $u$ by
$$u_*(x):=\inf\{t\geq 0 \,:\, \hbox{meas}(U^t)<\omega_n|x|^n\},$$
where $U^t$ are the level sets of $u=|u|$, that is, $U^t=\{x\in\mathbb{R}^n\,:\, |u(x)|>t\}$, and
$\omega_n$ denotes the   measure  of the standard unit  ball of $\rn$.
In particular, $(|x|^{-\alpha})_*=|x|^{-\alpha}$ for all $\alpha>0$, see \cite[3.3--(ii)]{ll}.
By the well known  Polya--Szego inequality (see \cite{talenti} and \cite{pz})
$$\int_{\mathbb{R}^n}|\nabla u_*|^p dx\leq \int_{\mathbb{R}^n}|\nabla u|^p dx,$$
and $u_*\in D_1^p(\rn)$, being $\int_{\mathbb{R}^n}|  u_*|^{\crit} dx=\int_{\mathbb{R}^n}| u|^{\crit} dx$.
Furthermore, by Theorem~3.4. of \cite{ll}
$$\int_{\rn}\frac{|u|^{p^\star(s)}}{|x|^s} \, dx\leq \int_{\rn}\frac{|u_*|^{p^\star(s)}}{|x|^s} \, dx\quad
\hbox{ and}\quad\int_{\rn}\frac{|u|^p}{|x|^p}\,  dx\leq \int_{\rn}\frac{|u_*|^p}{|x|^p}\,  dx.$$
Combining the above inequalities and the fact that $\mu\geq 0$,
we get that also $u_*$ is a minimizer and achieves the infimum of \eqref{def:Ks}.
Hence the equality sign   holds in all the inequalities above. In particular,
$$\int_{\rn}\frac{|u|^{p^\star(s)}}{|x|^s} \, dx=
\int_{\rn}\frac{|u_*|^{p^\star(s)}}{|x|^s} \, dx\quad
\hbox{ and}\quad\mu\int_{\rn}\frac{|u|^p}{|x|^p}\,  dx=\mu \int_{\rn}\frac{|u_*|^p}{|x|^p}\,  dx.$$
From Theorem~3.4 of \cite{ll}, in the case of equality, it then follows that $u=|u|=u_*$
if either $\mu\neq 0$ or if $s\neq 0$.
In particular, $u$ is positive, radially symmetric and decreasing with respect to 0. Hence $u$ must approach
a limit as $|x|\to\infty$, which must be zero, being $u\in L^{\crit}(\rn)$.
\hfill$\Box$

\section{Appendix 2: The case $\mu<0$}\label{sec:neg:mu}
As mentioned above, when $s=0$ and $\mu<0$, there is no extremal for \eqref{def:Ks}. More precisely, we have the following:

\medskip\noindent{\bf Claim \ref{sec:neg:mu}.1:} {\it Condition $\mu\leq 0$ entails that
$$K(n,p,\mu,0)=K(n,p,0,0).$$
In particular, there are no extremals when $\mu<0$.}\par
\smallskip\noindent{\it Proof of Claim \ref{sec:neg:mu}.1:} Since $\mu\leq 0$, we have that
\begin{equation}\label{ineq:1:neg:mu}
K(n,p,\mu,0)^{-1}\geq K(n,p,0,0)^{-1}.
\end{equation}
Let $u\in\dunp\setminus\{0\}$ be an extremal for $K(n,p,0,0)^{-1}$.
Fix $\alpha\in\rr$ and let $e_1$ be a nontrivial vector of $\rn$. We
define
\begin{equation}
\label{def:ualpha:neg:mu} u_\alpha(x):=u(x-\alpha e_1)
\end{equation}
for all $x\in\rn$. With a change of variables, we have
$$\frac{\int_{\rn}|\nabla u_\alpha|^p dx-\mu\int_{\rn}{|u_\alpha|^p}|x|^{-p} dx}
{\left(\int_{\rn}|u_\alpha|^{\crit}\right)^{{p}/{\crit}}}=
\frac{\int_{\rn}|\nabla u|^p dx-\mu\int_{\rn}{|u|^p}|x+\alpha
e_1|^{-p} dx} {\left(\int_{\rn}|u|^{\crit}\right)^{{p}/{\crit}}},$$
so that
$$\lim_{\alpha\to \infty}\frac{\int_{\rn}|\nabla u_\alpha|^p dx-
\mu\int_{\rn}{|u_\alpha|^p}|x|^{-p}
dx}{\left(\int_{\rn}|u_\alpha|^{\crit}\right)^{{p}/{\crit}}}=
\frac{\int_{\rn}|\nabla u|^p
dx}{\left(\int_{\rn}|u|^{\crit}\right)^{{p}/{\crit}}}=\frac1{K(n,p,0,0)}.$$
Therefore, $K(n,p,\mu,0)^{-1}\leq K(n,p,0,0)^{-1}$. Combining this with \eqref{ineq:1:neg:mu}, we obtain that $K(n,p,\mu,0)^{-1} =
K(n,p,0,0)^{-1}$. This proves Claim~\ref{sec:neg:mu}.1. \hfill$\Box$

\medskip\noindent Taking $u$ an extremal for $K(n,p,0,0)^{-1}$ and $u_\alpha$ as in \eqref{def:ualpha:neg:mu}, we get after some computations that
$$\max_{t\geq 0}\Phi(t u_\alpha)<\frac{1}{n}K(n,p,\mu,0)^{-n/p}$$
for $\alpha$ large when $0< s<\min\{p, (n-p)/(p-1)\}$. {\it This
permits to extend the proof given in Sections~$2$ and~$3$ to the
case $\mu< 0$ and $0<s<\min\{p, (n-p)/(p-1)\}$}.

\medskip\noindent We present here an alternative approach that allows to recover the full range $\mu<\mu_1$. Define
$$D_{1,r}^p(\rn):=\{u\in D_1^p(\rn)/\, u\hbox{ is radially symmetrical}\}$$
and for all $p\in (1,n)$, $s\in (0,p)$ and $\mu<\mu_1$, we let
\begin{equation}\label{def:Ks:r}
\frac{1}{K_r(n,p,\mu,s)}:=\inf_{u\in
D_{1,r}^p(\rn)\setminus\{0\}}\frac{\int_{\rn}|\nabla u|^p\, dx
-\mu\int_{\rn}{|u|^p}|x|^{-p}\,
dx}{\left(\int_{\rn}{|u|^{\crits}}|x|^{-s}\,
dx\right)^{{p}/{\crits}}}.
\end{equation}
Arguing as in Section \ref{sec:app}, we have
\begin{prop} For all $p\in (1,n)$, $s\in (0,p)$ and $\mu<\mu_1$, there are nonnegative extremals for $K_r(n,p,\mu,s)^{-1}$.
\end{prop}
In particular, a consequence of Theorem \ref{th:ext} and the remarks following this theorem is that
$$\begin{gathered}
K(n,p,\mu,s)=K_r(n,p,\mu,s),\\
\hbox{ when  }\mu\in [0,\mu_1)\hbox{ and }s\in[0,p),\hbox{ with }\mu+s>0;
\end{gathered}$$
while
$$K(n,p,\mu,s)>K_r(n,p,\mu,s)\hbox{ when }\mu<0\hbox{ and }s\in (0,s_\mu).$$

\medskip\noindent Since we have the existence of extremals in the radial case, one
can carry out the proofs of Sections~2 and~3 by restricting to radial functions and by replacing $K(n,p,\mu,s)$
in the definition  \eqref{bnd:c} of $c_*$ by $K_r(n,p,\mu,s)$.
This proves Theorem~\ref{th:main} in the case $\mu<0$.

\bigskip\noindent{\it Acknowledgements.}
  This work has been started while
F.~Robert was visiting the {\it Universit\`a degli Studi} of Perugia
in May 2007 with a GNAMPA--INdAM visiting professor position. He
thanks the members of the Mathematics department of the University
of Perugia for their hospitality and GNAMPA for the support. He is
also supported by Grants 2-CEx06-11-18/2006 and CNCSIS A--589/2007.
The first two authors were supported by the Italian MIUR project
titled {\it ``Metodi Va\-ria\-zio\-na\-li ed Equazioni Differenziali
non Lineari''}.\par
\noindent The authors thank the referee for very helpful and constructive comments on the manuscript.


\begin{thebibliography}{12}

\bibitem{afp} {\sc B. Abdellaoui, V. Felli and I. Peral}, Existence and non--existence
 results for quasilinear elliptic equations involving the
 $p$--Laplacian, {\it Boll. UMI}, {\bf 9--B} (2006), 445--484.

\bibitem{ar}
{\sc A. Ambrosetti and  P.H. Rabinowitz}, Dual variational methods
in critical point theory and applications, {\it J.  Funct. Anal.}, {\bf 14} (1973), 349--381.

\bibitem{aubin}
{\sc T. Aubin}, Probl\`emes isop\'erim\'etriques et espaces de
Sobolev, {\it J. Differential Geom.}, {\bf 11} (1976), 573--598.

\bibitem{ch} {\sc D. Cao and P. Han}, Solutions for semilinear
elliptic equations with critical exponents and Hardy potential, {\it
J. Differ. Equations}, {\bf205} (2004), 521--537.

\bibitem{ck} {\sc D. Cao and D. Kang}, Solutions of quasilinear elliptic problems involving a
Sobolev exponent and multiple Hardy--type terms, {\it J. Math. Anal.
Appl.}, {\bf 333} (2007), 889--903.


\bibitem{cw}
{\sc F. Catrina and Z.Q. Wang}. On the Caffarelli--Kohn--Nirenberg
inequalities: sharp constants, existence (and nonexistence), and
simmetry of extremal functions, {\it Comm. Pure Appl. Math.}, {\bf 54}
(2001), 229--258.

\bibitem{dh}
{\sc F. Demengel and E. Hebey}, On some nonlinear equations
involving the $p$--Laplacian with critical Sobolev growth, {\it Adv.
Differential Equations}, {\bf 3} (1998),  533--574.

\bibitem{druet}
{\sc O. Druet}, Generalized scalar curvature type equations on
compact Riemannian manifolds, {\it Proc. Roy. Soc. Edinburgh Sect.
A}, {\bf 130} (2000),   767--788.

\bibitem{evans}
{\sc L. C. Evans}, Weak convergence methods for nonlinear partial
differential equations, {\it CBMS Regional Conference Series in
Mathematics}, {\bf 74}. Published for the Conference Board of the
Mathematical Sciences, Washington, DC; by the American Mathematical
Society, Providence, RI, 1990. viii+80 pp.

\bibitem{ft} {\sc V. Felli and S. Terracini},
Elliptic equations with multi--singular inverse--square potentials
and critical nonlinearity, {\it  Comm. Partial Differential
Equations}, {\bf 31} (2006), 469--495.

\bibitem{gm} {\sc M. Gazzini and R. Musina}, On a Sobolev type inequality related to the weighted $p$--Laplace
operator, to appear in {\it J. Math. Anal. Appl.}


\bibitem{gr1} {\sc N. Ghoussoub and F.Robert}, The effect of curvature on the best constant in the Hardy--Sobolev inequality,
{\it Geom. Funct. Anal.}, {\bf 16} (2006), 897--908.


\bibitem{gr2} {\sc N. Ghoussoub and F.Robert}, Concentration estimates for Emden--Fowler equations
with boundary singularities and critical growth,
{\it  Int. Math. Res. Pap. IMRP}, Vol. 2006,
Article ID 21867, 1--85.


\bibitem{gy}
{\sc N. Ghoussoub and C. Yuan}, Multiple solutions for quasi--linear
PDEs involving the critical Sobolev and Hardy exponents, {\it Trans.
Amer. Math. Soc.}, {\bf 352} (2000),   5703--5743.

\bibitem{gv}
{\sc M. Guedda and L. V\'eron}, Quasilinear elliptic equations
involving critical Sobolev exponents, {\it  Nonlinear Anal.}, {\bf
13} (1989),   879--902.

\bibitem{h} {\sc P. Han}, Quasilinear elliptic problems with critical exponents and Hardy
terms, {\it Nonlinear Anal.}, {\bf 61} (2005), 735--758.


\bibitem{kl}
{\sc D. Kang and G. Li}, On the elliptic problems involving multi--singular
inverse square potentials and multi--critical Sobolev--Hardy exponents, {\it  Nonlinear Anal.}, {\bf 66} (2007),
1806--1816.




\bibitem{kp0}
{\sc D. Kang and S. Peng}, Existence of solutions for elliptic
problems with critical Sobolev--Hardy exponents, {\it Israel J.
Math.}, {\bf 143} (2004), 281--297.


\bibitem{kp}
{\sc D. Kang and S. Peng}, Solutions for semilinear elliptic
problems with critical Sobolev--Hardy exponents and Hardy potential,
{\it Appl. Math. Lett.}, {\bf 18} (2005), 1094--1100.

\bibitem{li} {\sc J. Li}, Equation with critical Sobolev--Hardy exponents, {\it Intern. J. Math.
Math. Sci.}, {\bf20} (2005), 3213--3223.

\bibitem{ll} {\sc E.H. Lieb and M. Loss},  Analysis, {\it Graduate
Studies in Mathematics}, {\bf14}, Amer. Math. Soc., Providence, RI, 1997. xviii+278 pp.

\bibitem{lions}
{\sc P. L. Lions}, The concentration--compactness principle in the
calculus of variations. The limit case. I, II, {\it Rev. Mat.
Iberoamericana}, {\bf 1} (1985), 145--201 and 45--121.

\bibitem{pz} {\sc G. P\'olya and G. Szeg\"o},
Isoperimetric Inequalities in Mathematical Physics, {\it Annals of Mathematics Studies},
{\bf27}, Princeton University Press, Princeton, N. J., 1951. xvi+279 pp.

\bibitem{patraf} {\sc P. Pucci and R. Servadei}, Existence, non--existence and
regularity of radial ground states for $p$--Laplacian equations with
singular weights,  {\em Ann. Inst. H. Poincar\'e
A.N.L.}, {\bf 25} (2008), 505--537.


\bibitem{pucciservadei} {\sc P. Pucci and R. Servadei}, Regularity
of Weak Solutions of Homogeneous or Inhomogeneous Quasilinear
Elliptic Equations, to appear in {\em  Indiana Univ. Math. J}.

\bibitem{rodemich}
{\sc E. Rodemich}, The Sobolev inequalities with best possible
constants, {\it  Analysis Seminar at California Institute of
Technology} (1966).

\bibitem{saintier1}
{\sc N. Saintier}, Asymptotic estimates and blow--up theory for
critical equations involving the $p$--Laplacian, {\it Calc. Var.
Partial Differential Equations}, {\bf 25} (2006),   299--331.

\bibitem{saintier2}
{\sc N. Saintier}, Blow--up theory for symmetric critical equations
involving the $p$--Laplacian. {\it  NoDEA, Nonlinear Differ.
Equ. Appl.}, {\bf15} (2008), 227--245.

\bibitem{serrin}
{\sc J. Serrin}, Local behavior of solutions of quasi--linear
equations, {\it  Acta Math.}, {\bf 111} (1964), 247--302.

\bibitem{struwe}
{\sc M. Struwe}, Variational methods. Applications to nonlinear
partial differential equations and Hamiltonian systems, Third
edition. {\it Ergebnisse der Mathematik und ihrer Grenzgebiete.} 3.
Folge. A Series of Modern Surveys in Mathematics, {\bf 34}.
Springer--Verlag, Berlin, 2000. xviii+274 pp.

\bibitem{talenti}
{\sc G. Talenti}, Best constant in Sobolev inequality, {\it  Ann.
Mat. Pura Appl.}, {\bf 110} (1976), 353--372.

\bibitem{terr}
{\sc S. Terracini}, On positive solutions to a class equations with
a singular coefficient and critical
 exponent, {\it Adv. Differential Equations}, {\bf 2} (1996), 241--264.

\bibitem{tolk}
{\sc P. Tolksdorf}, Regularity for a more general class of
quasilinear elliptic equations, {\it J. Differential Equations},
{\bf 51} (1984),   126--150.

\bibitem{vaz}
{\sc J. L. V\`azquez}, A strong maximum principle for some
quasilinear elliptic equations, {\it Appl. Math. Optim.}, {\bf 12}
(1984), 191--202.
\end{thebibliography}
\end{document}